\documentclass[10pt,twoside,reqno]{amsart} %%%%%%%%%%%%%%%%%% reqno: c'est right equation number
\usepackage{amsmath, amsfonts, amsthm, amssymb}
\usepackage[francais,english]{babel}
\usepackage[latin1]{inputenc}
\usepackage{amsmath}
\usepackage{amsfonts}
\usepackage{amsthm}
\usepackage[]{graphics}
\usepackage{color}
\usepackage[T1]{fontenc}
\usepackage{setspace}
\usepackage{alltt}
\usepackage{algorithm}
\usepackage{algorithmic}

\usepackage[T1]{fontenc}
\usepackage[latin1]{inputenc}
\usepackage{amsmath,latexsym}
\usepackage{here}
\usepackage{graphicx}
\usepackage{slashbox}
\usepackage{multirow}
\usepackage{amssymb}

\definecolor{red}{rgb}{0.82,0.00,0.00}  % rouge

\numberwithin{equation}{section}
\input epsf
 

\newcommand{\ds}{\displaystyle}
\def\R{{\rm I\hspace{-0.50ex}R} }
\def\N{{\rm I\hspace{-0.50ex}N} }

\def\0{{\bf 0}}

\def\E{{\mathrm{E}}}
\def\V{{\mathrm{Var}}}

\newtheorem{lem}{Lemma}[section]

\newtheorem{rmq}[lem]{Remark}

\renewcommand{\epsilon}{\varepsilon}
\def\twoplot[#1]#2#3#4#5{
\begin{figure}[hbt]
\begin{multicols}{2}
\begin{center}
    \includegraphics*[#1]{#2}
    \caption{\label{#2} #4}
\end{center}
\begin{center}
    \includegraphics*[#1]{#3}
    \caption{\label{#3} #5}
\end{center}
\end{multicols}
\end{figure}
}

\begin{document}
\bibliographystyle{plain}

\title[ Adaptive Monte Carlo method] {Adaptive stratified Monte Carlo algorithm for numerical computation of
integrals}

\author[Toni Sayah  ]{{Toni Sayah $^{\dagger}$} }
\thanks{ \today.
\newline
$^{\dagger}$  Unit\'e de recherche EGFEM, Facult\'e des sciences,
Universit\'e Saint-Joseph, Lebanon.
\newline
toni.sayah@usj.edu.lb.}

\maketitle

\begin{abstract}

In this paper, we aim to compute numerical approximation integral
by using an adaptive Monte Carlo algorithm. We propose a
stratified sampling algorithm  based on an iterative method which
splits the strata following some quantities called indicators
which indicate where the variance takes relative big values. The
stratification method is based on the optimal allocation strategy
in order to decrease the variance from iteration to another.
Numerical experiments show and confirm the efficiency of our
algorithm.

\mbox{}

\end{abstract}

\section{Introduction}\label{intro}
This paper deals with adaptive Monte Carlo method (AMC) to
approximate the integral of a given function $f$ on the hypercube
$[0,1[^d, d\in \N^*$. The main idea is to guide the random points
in the domain in order to decrease the variance and to get better
results. The corresponding algorithm couples two methods: the
optimal allocation strategy and the adaptive stratified sampling.
In fact, it proposes to split the domain into separate regions
(called mesh) and to use an iterative algorithm which calculates
the number of samples in every region by using the optimal
allocation strategy and then refines the parts of the mesh
following some quantities called indicators which indicate where
the variance takes a relative big values. \\

A usual technique for reducing the mean squared error of a
Monte-Carlo estimate is the so-called stratified Monte Carlo
sampling, which considers sampling into a set of strata, or
regions of the domain, that form a partition (a stratification) of
the domain (see \cite{PG} and the references therein for a
presentation more detailed). It is efficient to stratify the
domain, since when allocating to each stratum a number of samples
proportional to its measure, the mean squared error of the
resulting estimate is always smaller or equal to the one of the
crude Monte-Carlo estimate. For a given partition of the domain
and a fixed total number of random points, the choice of the
number of samples in each stratum is very important for the
results and precision. The { optimal allocation strategy} (see for
instance \cite{PEBJ} or \cite{ACRM}) allows to get the better
distribution of the samples in the set of strata in order to
minimize the variance. We give in the next section a brief summary
of the this strategy which will
be the basic tools of our adaptive algorithm.\\

In the other hand, it is important to stratify the domain in
connection with the function $f$ to be integrated and to allocate
more strata in the region where $f$ has larger local variations.
Many research works propose multiple methods and technics to
stratify the domain: \cite{ACRM} for the adaptive stratified
sampling for Monte-Carlo integration of differentiable functions,
\cite{CLSM} for the adaptive integration and approximation over
hyper-rectangular regions with applications to basket option
pricing, \cite{PFBT}, \dots . \\

The paper is organized as follows. Section 2 describes the
adaptive method. We begin by giving a summarize of the { optimal
allocation strategy} and then describe the adaptive algorithm
consisting in stratifying the domain. In section 3, we perform
numerical investigations showing the powerful of the proposed
adaptive algorithm.
\section{Description of the adaptive algorithm}
\noindent In this section, we will describe the AMC algorithm
which is based on indicators to guide the repartition of the
random points in the domain. In our algorithm, the indicators are
based on an approximation of the variance expressed on different
regions in the domain. We detect those where the indicators are
bigger than their mean value up to a constant and we
split them in small regions.\\

\subsection{Optimal choice of the numbers of samples}
\noindent
Let    $D=[0,1[^d$ be the unit hypercube of $\R^d$, $d\in\N^*$,  and  $f:D\rightarrow \R$ a Lebesgue-integrable    function.
We want to estimate
\[
\ds \mathcal{I}(f) = \int_D f(x) d \lambda(x),
\]
where $\lambda$ is the Lebesgue measure on $\R^d$.

The classical MC estimator  of $\mathcal{I}(f)$ is
\[
\ds \mathcal{\bar{I}}_{\text{MC}}(f) = \ds \frac{1}{N}
\sum_{i=1}^{N} f \circ U_i,
\]
where  $U_i, 1\le i\le N$,  are
independent random variables uniformly distributed over $D$.  $\mathcal{\bar{I}}_{\text{MC}}(f)$ is an unbiased estimator of $\mathcal{I}(f)$,
which  means that $\E[\mathcal{\bar{I}}_{\text{MC}}(f)]=\mathcal{I}(f)$.
Moreover, if $f$ is square-integrable,  the variance of $\mathcal{\bar{I}}_{\text{MC}}(f)$ is
\[
\V(\mathcal{\bar{I}}_{\text{MC}}(f) )=
\ds \frac{\sigma^2(f)}{N}
\]
where
\[
\sigma^2(f) = \ds \int_{D} f(x)^2 d\lambda(x) - \left(\ds \int_{D}
f(x) d\lambda(x) \right)^2.
\]
Variance reduction techniques aim to produce alternative
estimators having smaller variance than crude MC. Among these
techniques, we focus on stratification strategy.
The idea is to split   $D$   into separate regions, take a sample
of points from each such region, and combine the results to
estimate $\mathcal{I}(f)$. Let $\{D_1,\ldots,D_p\}$ be a partition
of $D$. That is a set of sub-domains such that
\[
D=\bigcup_{i=1}^p D_i \quad \text{and}\quad D_i\cap D_j =\emptyset \ \text{for}\ i\neq j.
\]
We consider $p$ corresponding integers $n_1,\ldots, n_p$.
Here, $n_i$ will be the number of samples to draw from $D_i$.
For $1\leq i\leq p$, let  $a_i = \ds \int_{D_i} d \lambda(x)$ be the measure of
 $D_i$ and $\mathcal{I}_{i}(f) = \ds \int_{D_i} f(x) d\lambda (x)$ be the
integral  of $f$ over $D_i$.
% $N= \ds \sum_{i=1}^p n_i$.
We have $\lambda(D) = \ds \sum_{i=1}^p a_i$ and  $ \mathcal{I}(f)=
\ds \sum_{i=1}^p \mathcal{I}_{i}(f)$. Furthermore, for   $1\le i
\le p$, let $\pi_i = \ds \frac{1_{D_i}}{a_i} \lambda$  be the
density function of the uniform distribution over $D_i$ and
consider a set of $n_i$ random variables $X^{(i)}_1,
\dots,X^{(i)}_{n_i}$ drawn from  $\pi_i$. We suppose that the
random variables  $X^{(i)}_j,1\leq j\leq n_i, 1\leq i\leq p$, are
mutually independent.

For $1\leq i\leq p$, let $S_i$ be the MC estimator of $\mathcal{I}_{i}(f) $
defined by:
\[
\ds S_i = \ds \frac{1}{n_i} \sum_{k=1}^{n_i} f \circ X^{(i)}_k.
\]
Then, the integral $\mathcal{I}(f)$ can be estimated by:
\[
 \mathcal{\bar{I}}_{\text{SMC}}(f)  = \ds \sum_{i=1}^p   a_i S_i = \ds \sum_{i=1}^p \frac{a_i}{n_i} \sum_{k=1}^{n_i} f \circ X^{(i)}_k.
\]
We call $\mathcal{\bar{I}}_{\text{SMC}}(f)$ the \emph{stratified Monte Carlo
 estimator} of $\mathcal{I}(f)$.
It is easy to show that $ \mathcal{\bar{I}}_{\text{SMC}}(f)$ is an unbiased estimator of $\mathcal{I}(f)$ and, if $f$ is square-integrable,  the variance of $\mathcal{\bar{I}}_{\text{SMC}}(f)$ is
\[
\V(\mathcal{\bar{I}}_{\text{SMC}}(f)) = \ds \sum_{i=1}^{n_i} a_i^2  \ds \frac{\sigma^2_i(f)}{n_i}
\]
where
\[
\sigma^2_i(f) = \ds \int_{D} f(x)^2 d\pi_i(x) - \left( \int_{D} f(x) d\pi_i(x) \right)^2,  \quad  \forall 1\leq i\leq p.
\]

The choice of the integers $n_i$, $i=1,\ldots,p$ is crucial in
order to reduce   $\V(\mathcal{\bar{I}}_{\text{SMC}}(f))$.  A
frequently made choice is proportional allocation which takes the
number $n_i$ of points in each sub-domain $D_i$ proportional to
its measure. In other words, if $N=\ds \sum_{i=1}^p n_i$, then
$n_i= Na_i, i=1,\ldots,p$. \\

\noindent For this choice, we have
\[
\V(\mathcal{\bar{I}}_{\text{MC}}(f))=\V(\mathcal{\bar{I}}_{\text{SMC}}(f))+
\frac{1}{N} \sum_{i=1}^p a_i \left(\frac{\mathcal{I}_i(f)}{a_i} - \mathcal{I}(f)\right)^2,
\]
 and hence, $\V(\mathcal{\bar{I}}_{\text{SMC}}(f))\leq\V(\mathcal{\bar{I}}_{\text{MC}}(f))$.

To get an even smaller variance, one can consider The {optimal
allocation} which aims  to minimize
\[
V(n_1,\ldots,n_p) = \ds \sum_{i=1}^{n_i} a_i^2  \ds
\frac{\sigma^2_i(f)}{n_i},
\]
as a function of $n_1,\ldots,n_p$, with  $N= \ds \sum_{i=1}^p
n_i$. Let
\[
\delta = \ds \frac{1}{N} \sum_{i=1}^p a_i \sigma_i(f).
\]
Using the inequality of Cauchy-Schwarz, we have
\begin{align*}
V\left(\frac{a_1 \sigma_1(f)}{\delta}, \dots , \frac{a_p \sigma_p(f)}{\delta}\right)
&= \frac{1}{N} \left(\sum_{i=1}^p a_i \sigma_i(f)\right)^2 \\
&\leq \frac{1}{N}\left(\sum_{i=1}^p \frac{a_i^2 \sigma_i(f)^2}{n_i}\right)
\sum_{i=1}^p n_i\\
& \leq V(n_1,\ldots, n_p).
\end{align*}
Hence, the optimal choice of $n_1,\ldots,n_p$ is given by
\begin{equation}\label{rept}
n_i= \frac{a_i\sigma_i(f)}{\delta}, \quad i=1, \dots,p.
\end{equation}

In order  to compute the number $n_i$ of random
points in    $D_i$ using (\ref{rept}),  one can approximate  $\sigma_i(f)$
 by:
\begin{equation}\label{approxsigmai}
\bar{\sigma}^2_i (f) =  \ds \frac{1}{n_i} \sum_{j=1}^{n_i} (f \circ X^{(i)}_j)^2 - \Bigl( \frac{1}{n_i} \sum_{j=1}^{n_i} f \circ X^{(i)}_j \Bigr)^2.
\end{equation}
 For $1\leq i\leq p$, we  will  denote
\begin{equation}\label{ni}
\bar{n}_i = \ds \frac{a_i \bar{\sigma}_i (f)}{\bar{\delta}}
\end{equation}
where
\begin{equation}\label{deltai}
\bar{\delta} = \ds \frac{1}{N} \sum_{i=1}^p a_i \bar{\sigma}_i
(f).
\end{equation}
\subsection{Description of the algorithm}
\noindent The adaptive MC  scheme  aims to guide, for a fixed
global number $N$ of random points in the domain $D$, the
generation of random points in every sub-domain in order to get
more precise estimation on the desired integration. It is based on
an iterative algorithm where the mesh (repartition of the
sub-domains in $D$) evolves with iterations. Let $L$ be the total
number of desired iterations and $D^\ell, 1\leq \ell\leq L,$ be
the corresponding mesh such that
\[
D^\ell=\bigcup_{i=1}^{p_\ell} D^\ell_i \quad \text{and}\quad
D^\ell_i\cap D^\ell_j =\emptyset \ \text{for}\ i\neq j,
\]
where $p_\ell$ is the number of the sub-domains in $D^\ell$. We
start the iterations with a subdivision of the domain $D^1=D$
using $p_1$ identical sub-cubes with given equal numbers of random
points $n_{i,1}$ in each sub-domain $D^1_i, i=1, \dots, n_{i,1}$,
such that $N=\ds \sum_{i=1}^{p_1} n_{i,1}$. \\
%The choice of the
%total number of iterations $L$ corresponds to the desired
%precision related to the calculated variance at the level $L$. \\

The main idea of the algorithm consists for every iteration $1\le
\ell \le L$, to refine some region $D^\ell_i$, $1\leq i\leq
p_\ell$, of the mesh $D^\ell$ where the function $f$ presents more
singularities (big values of the variance) and hence must be
better described. This technique is based on some quantities
called \emph{indicators} and denoted ${V}_{i,\ell}$  which give
informations about the contribution of $D^\ell_i$  in the
calculation of the variance of the MC method at this level,
approximated  by
\begin{equation}\label{numvariance}
{V}_\ell = \ds \sum_{i=1}^{p_l} {V}_{i,\ell}
\end{equation}
where
\begin{equation}\label{locnumvariance}
{V}_{i,\ell} = \ds \frac{a_{i,\ell}^2
\bar{\sigma}_{i,\ell}^2(f)}{\bar{n}_{i,\ell}}.
\end{equation}
Our goal is to decrease ${V}_\ell$ during the iterations. Then,
for every refinement iteration $\ell$ with a corresponding mesh
$D^{\ell}$ and corresponding numbers $\bar{n}_{i,\ell}$, we
calculate $\bar{\sigma}_{i,\ell}(f)$ and $\bar{\delta}_\ell$, and
update $\bar{n}_{i,\ell}$ by using the optimal choice of the
numbers of samples based on the formulas (\ref{approxsigmai}),
(\ref{deltai}) and (\ref{ni}) for all the sub-domains $D^\ell_i,
i=1, \dots, p_\ell$. For technical reason, we allow a minimal
number, denoted by $M_{rp}$ (practically we choose $M_{rp}=2$), of
random points in every sub-domain and then if $\bar{n}_{i,\ell} <
M_{rp}$ we set
$\bar{n}_{i,\ell} = M_{rp}$. \\
\noindent Next, we calculate the indicators ${V}_{i,\ell}$ and
${V}_\ell$, and then, we adapt the mesh $D^\ell$ to obtain the new
one $D^{\ell +1}$. The chosen strategy of the adaptive method
consists to mark the sub-domains $D^\ell_i$ such that
\[
{V}_{i,\ell} > C_m {V}^{mean}_\ell ,
\]
where $C_m$ is a positive constant bigger than $1$ and
${V}^{mean}_\ell$ is the mean value of ${V}_{i,\ell}$ defined as
\begin{equation}\label{meanval}
{V}^{mean}_\ell = \ds \frac{1}{p_\ell} \sum_{i=1}^{p_l}
V_{i,\ell},
\end{equation}
and to divide every marked sub-domains $D^\ell_i$ into small
parts, four equal sub-squares for $d=2$ and eight equal sub-cubes
for $d=3$, with equal number of random points in each part given
by
\[
\left\{
\begin{array}{rcl}
\medskip
\ds \max(\frac{\bar{n}_{i,\ell}}{4}, M_{pr}) \qquad \mbox{for
} d=2\\
\ds \max(\frac{\bar{n}_{i,\ell}}{8}, M_{pr}) \qquad \mbox{for }
d=3.
\end{array}
\right.
\]
{\rmq We stop the algorithm if the number of iterations reaches
$L$ or if the calculated variance is smaller that a tolerance
denoted by $\varepsilon$. We denote by $\ell_\varepsilon$ the
stopping iteration level of the following algorithm which
corresponds to a desired  tolerance $\varepsilon$  or at maximum
equals to $L$. }

\vspace{.3cm}
\noindent The algorithm can be  described as following : \\

\noindent {\underline {(Algo 1)}} : For a chosen $N$
with corresponding numbers $n_{i,1}$, and a given initial mesh
$D^1$ with corresponding sub-domains $D^1_i, i=1, \dots , p_1$,
{\small
\begin{alltt}
Generate \(n\sb{i,1},i=1,\dots,p\sb{1}\) random points \(X\sp{i}\sb{j},j=1,\dots,n\sb{i,1}\) in every sub-domain  \(D\sp{1}\sb{i}\).
set \(\ell=1\).
calculate \(V\sb{\ell}\) by using (\ref{numvariance}).
\textbf{While} \(l{\leq}L or V\sb{\ell}\le\varepsilon\)
    calculate \(\bar{\sigma}\sb{i,\ell}(f)\) and, \(\bar{\delta}\sb{\ell}\) and  update \(\bar{n}\sb{i,\ell},i=1,\dots,p\sb{\ell}\) by using (\ref{approxsigmai}), (\ref{deltai}) and (\ref{ni}).
    Generate corresponding random points \(X\sp{i}\sb{j},j=1,\dots,n\sb{i,l}\) in each sub-domain \(D\sp{\ell}\sb{i},i=1,\dots,p\sb{\ell}\).
    calculate \({V}\sb{i,\ell},i=1,\dots,p\sb{l}\) and \({V}\sp{mean}\sb{\ell}\) by using (\ref{locnumvariance}) and (\ref{meanval}).
    for (\(i=1:p\sb{\ell}\))
        if (\(V\sb{i,\ell}\ge C\sb{m} V\sp{mean}\sb{\ell}\))
            Divide the sub-domain \(D\sb{i}\sp{\ell}\) in \(m\) small parts (\(m=4\) in 2D and \(m=8\) in 3D).
            Associate to every one of this small parts the number of random points max(\(\ds\frac{\bar{n}\sb{i,l}}{m},M\sb{pr}\)).
            set \(p\sb{\ell}=p\sb{\ell}+m\).
        end if
    end for
    \(\ell=\ell+1\).
end loop
\(\ell\sb{\varepsilon}=\ell-1\).
calculate the adapt MC approximation \(\mathcal{I}\sb{AMC}=\ds\sum\sb{i=1}\sp{p\sb{\ell\sb{\varepsilon}}}\frac{a\sb{i,\ell\sb{\varepsilon}}}{\bar{n}\sb{i,\ell\sb{\varepsilon}}}\sum\sb{k=1}\sp{n\sb{i,\ell\sb{\varepsilon}}}f\circ{X\sp{i}\sb{k}}\).
\end{alltt}
}
The previous algorithm calculate an approximation of
${\mathcal{I}}(f)$ with an adaptive Monte Carlo method. If we are
interested by the numerical variance, we repeat the previous
algorithm $N_{ess}$ times and approximate the ${\mathcal{I}}(f)$
by the corresponding mean value
\[
\bar{\mathcal{I}}_{AMC} = \ds \frac{1}{N_{ess}}\sum_{i=1}^{N_{ess}} {\mathcal{I}}^i_{AMC},
\]
where ${\mathcal{I}}^i_{AMC}$ corresponds to the $i^{th}$ essay using (Algo 1). \\
The estimated variance will by given by the formula
\[
V_{AMC} = \ds \frac{1}{N_{ess} - 1} \Big(
\sum_{i=1}^{N_{ess}} ({\mathcal{I}}^i_{AMC})^2 - N_{ess} \, \bar{\mathcal{I}}_{AMC}^2
\Big).
\]
In fact, it is useless to repeat the (Algo 1) $N_{ess}$ times to
calculate $\bar{\mathcal{I}}_{AMC}$ and $V_{AMC}$, and it is
expensive for the CPU time. We can reduce the coast by running
(Algo 1) one time to define the mesh and to get
${\mathcal{I}}^1_{AMC}$ and then, we use the corresponding
sub-domains $D^{\ell_\varepsilon}_i, i=1, \dots ,
p_{\ell_\varepsilon}$ with the corresponding number of random
points $n_{i,\ell_\varepsilon},i=1,\dots, \ell_\varepsilon$ to
perform the rest of calculations ($N_{ess} -1$ essays).  The
corresponding algorithm can be describe as follow : \\

\noindent {\underline {(Algo 2)}} :
{\small
\begin{alltt}
Call algorithm (Algo 1) to define the mesh \(D\sp{\ell\sb{\varepsilon}}\sb{i},i=1,\dots,p\sb{\ell\sb{\varepsilon}}\) and calculate \(\mathcal{I}\sp{1}\sb{AMC}\)
Set \(n\sb{e} = 2\)
\textbf{While}  \(n\sb{e}{\leq}N\sb{ess}\)
    Generate corresponding random points \(X\sp{i}\sb{j},j=1,\dots,n\sb{i,\ell\sb{\varepsilon}}\) in each sub-domain \(D\sp{\ell\sb{\varepsilon}}\sb{i},i=1,\dots,p\sb{\ell\sb{\varepsilon}}\).
    Calculate \(\mathcal{I}\sp{n\sb{e}}\sb{AMC}=\ds\sum\sb{i=1}\sp{p\sb{\ell\sb{\varepsilon}}}\frac{a\sb{i,\ell\sb{\varepsilon}}}{\bar{n}\sb{i,\ell\sb{\varepsilon}}}\sum\sb{k=1}\sp{n\sb{i,\ell\sb{\varepsilon}}}f\circ{X\sp{i}\sb{k}}\)
    Set \(n\sb{e}=n\sb{e}+1\)
end loop calculate
\(\bar{\mathcal{I}}\sb{AMC}=\ds\frac{1}{N\sb{ess}}\sum\sb{i=1}\sp{N\sb{ess}}\mathcal{I}\sp{i}\sb{AMC}\) calculate
\(V\sb{AMC}=\ds\frac{1}{N\sb{ess}-1}\Big(\sum\sb{i=1}\sp{N\sb{ess}}({\mathcal{I}}\sp{i}\sb{AMC})\sp{2}-N\sb{ess}\,\bar{\mathcal{I}}\sb{AMC}\sp2\Big)\)
\end{alltt}
}
\section{Numerical experiments}
\noindent In this section, we perform in MATLAB several numerical
experiments to validate our approach and we compare between the MC and AMC methods. \\
\subsection{2D validations}
\begin{figure}[h]%\vspace{-0cm}
\includegraphics[height=3.5cm,width=4cm]{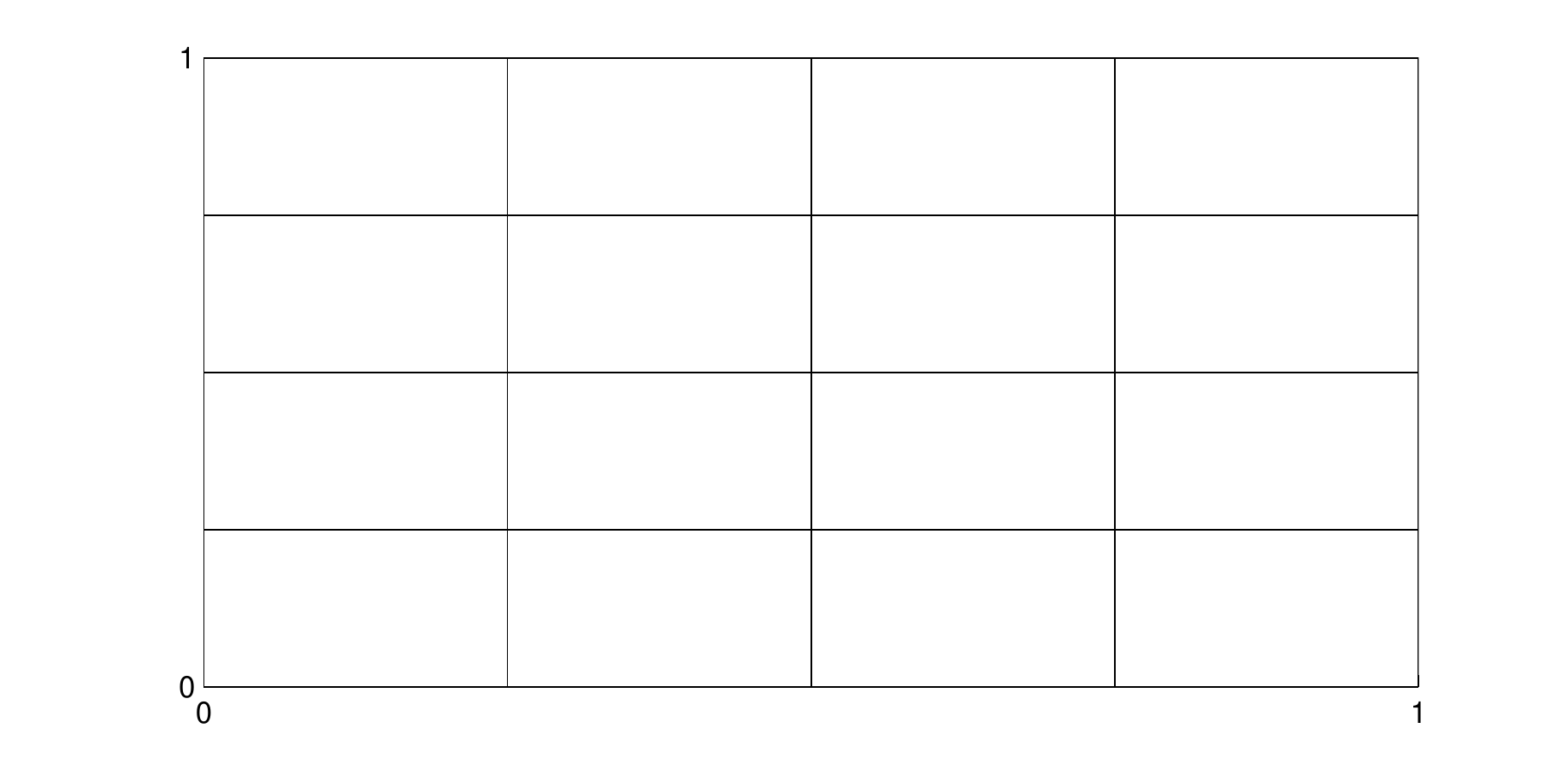}\vspace{-.6cm}
\caption{Initial partition $D^1_i, i=1,\dots, p_1$ ($p_1=16$)with
$N_0=4$.}\label{initpart}
\end{figure}

\noindent We consider the unit square $D=[0,1[^2$, $C_m=2$,
$M_{pr}=2$ and $\ell_\varepsilon=L$. The initial mesh is
constituted by a regular partition with $N_0=4$ segments in every
side of $D^1=D$ (see figure \ref{initpart}).

In this section, we show two particular cases of the function $f$.
The first treats an integrable but not continuous function which
presents a discontinuity along the border of the unit disc. The
second one treats a function concentrated in a part of $D$ and
vanishes in the rest on this domain. Both examples show the
powerful of the proposed AMC method.
\subsubsection{First test case} $\mbox{ }$\\
\noindent For the first test case, we consider  the function $f_c$
given on $D$ by
\[
f_c(x,y) = \left\{
\begin{array}{rcl}
1 \quad && \mbox{ if } x^2 + y^2 \le 1 \\
0 \quad && \mbox{ elsewhere. }
\end{array}
\right.
\]
The exact integration of $f_c$ over $D$ is equal to
$$I=\ds \int_D f_c(x,y) dx dy = \ds \frac{\pi}{4},$$
which is the quarter of the surface of the unit disc.\\

\begin{figure}[h!]

 \begin{minipage}[h]{.46\linewidth}\vspace{-.6cm}
  \centering\includegraphics[width=7.5cm]{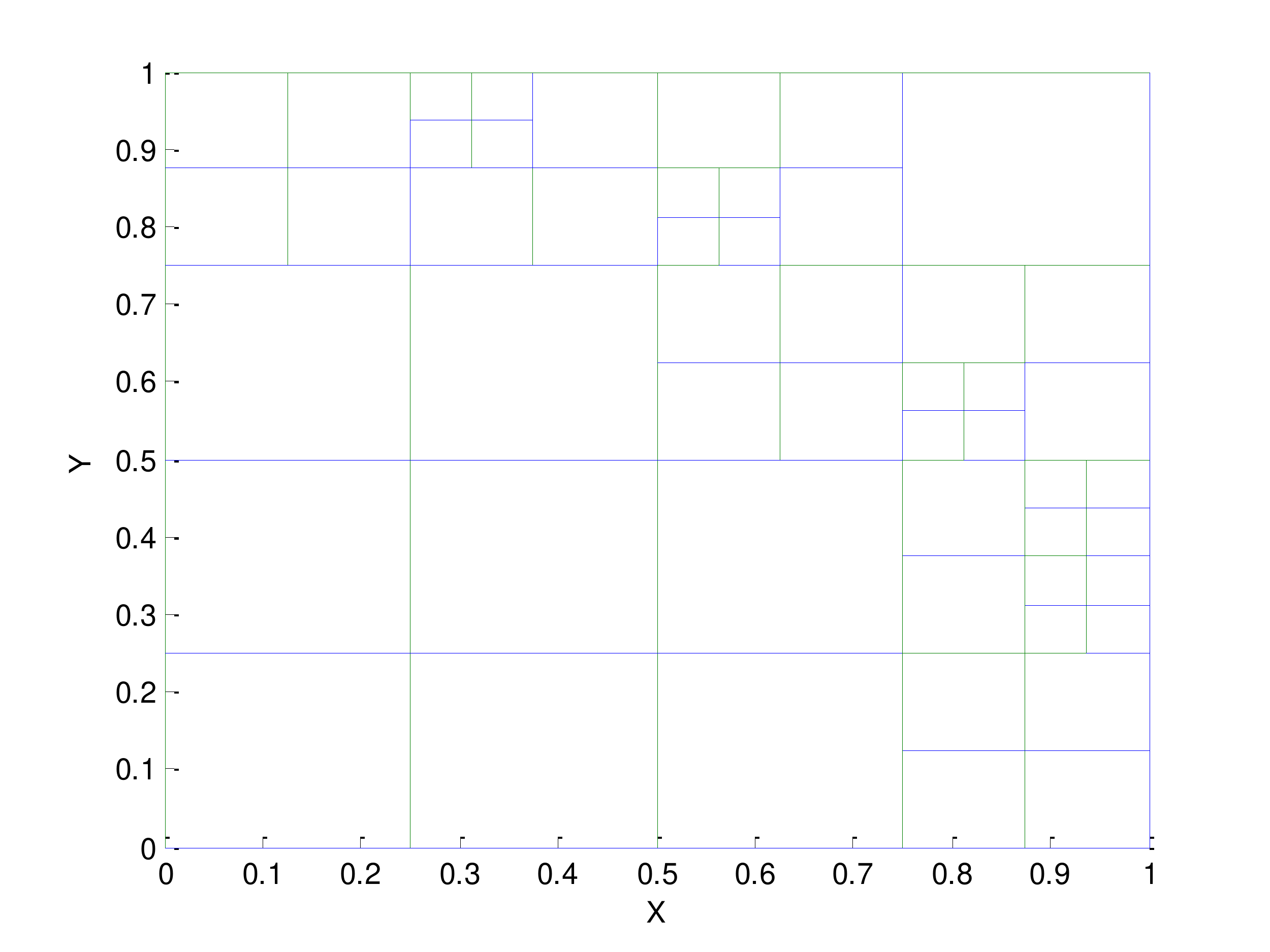}\vspace{-.7cm}
  \caption{Mesh for the second iteration.}  \label{step1}
 \end{minipage} \hfill
 \begin{minipage}[h]{.46\linewidth} \vspace{-.3cm}
  \centering\includegraphics[width=7.5cm]{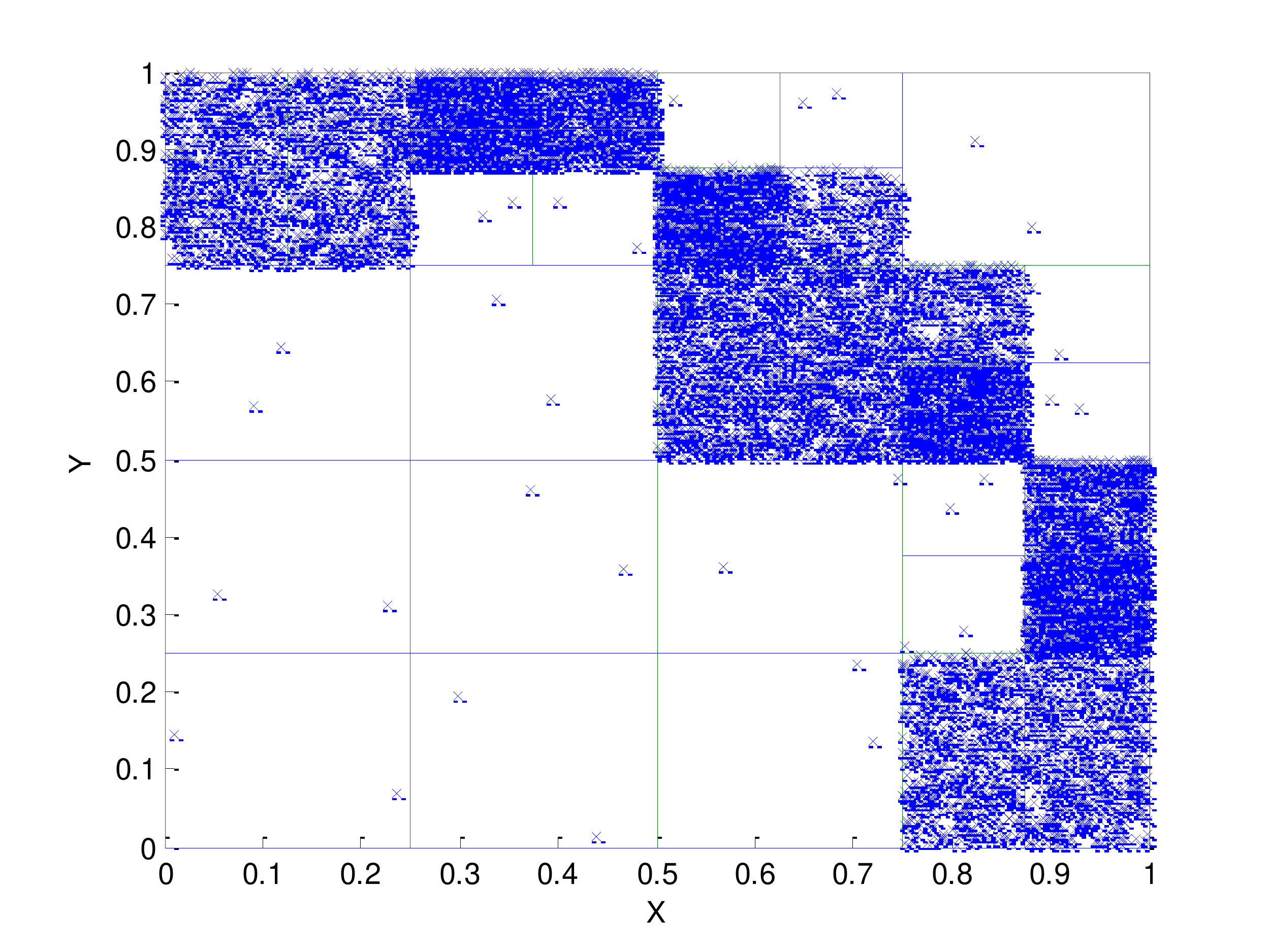}\vspace{-.7cm}
  \caption{repartition of the points  (second iteration).}
 \end{minipage}

\end{figure}

\begin{figure}[h!]
%\vspace{.5cm}
 \begin{minipage}[h]{.46\linewidth}%\vspace{-.9cm}
  \centering\includegraphics[width=7.5cm]{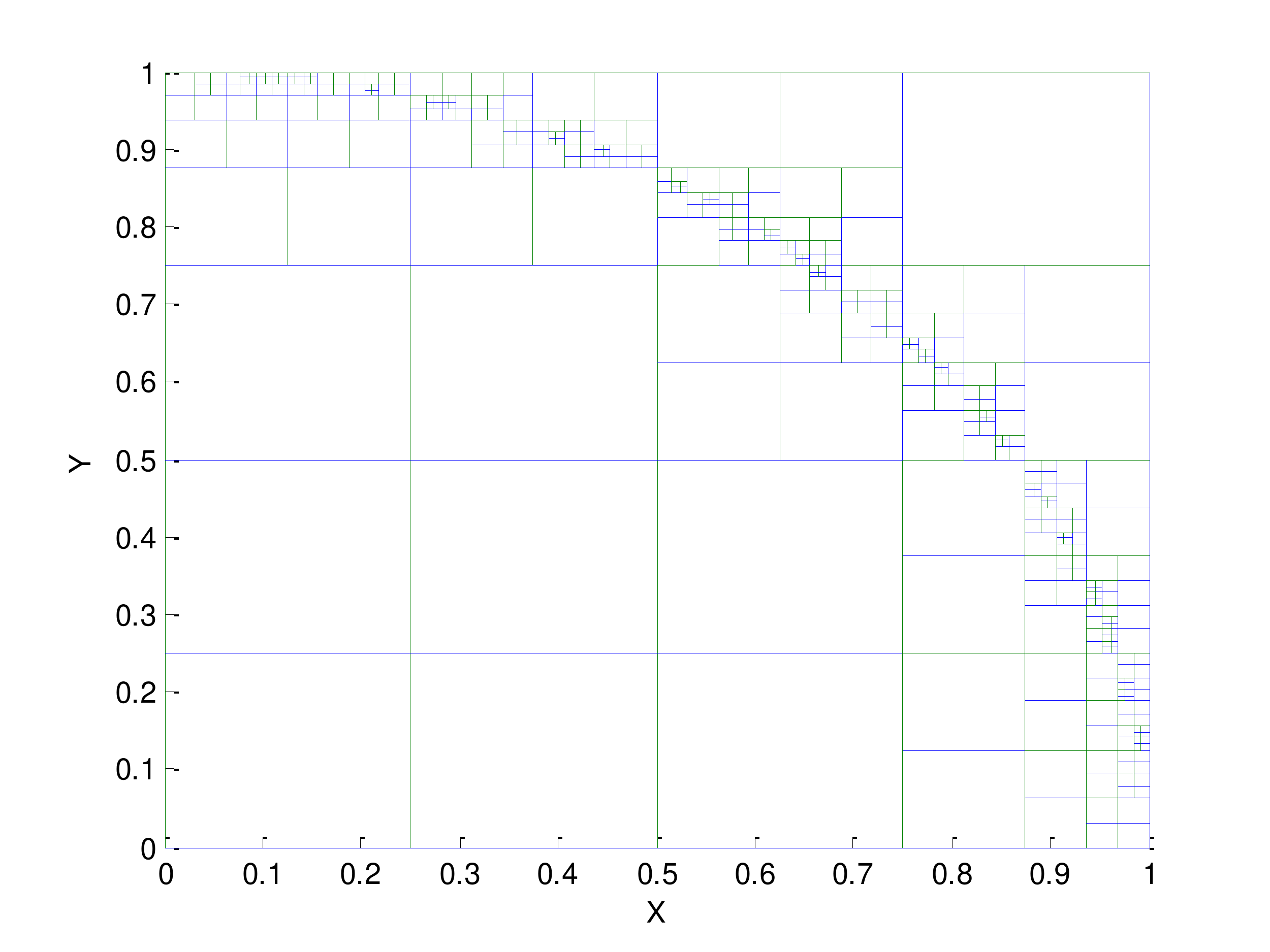}\vspace{-.7cm}
  \caption{Mesh for the sixth iteration.}
 \end{minipage} \hfill
 \begin{minipage}[h]{.46\linewidth}\vspace{.2cm}
  \centering\includegraphics[width=7.5cm]{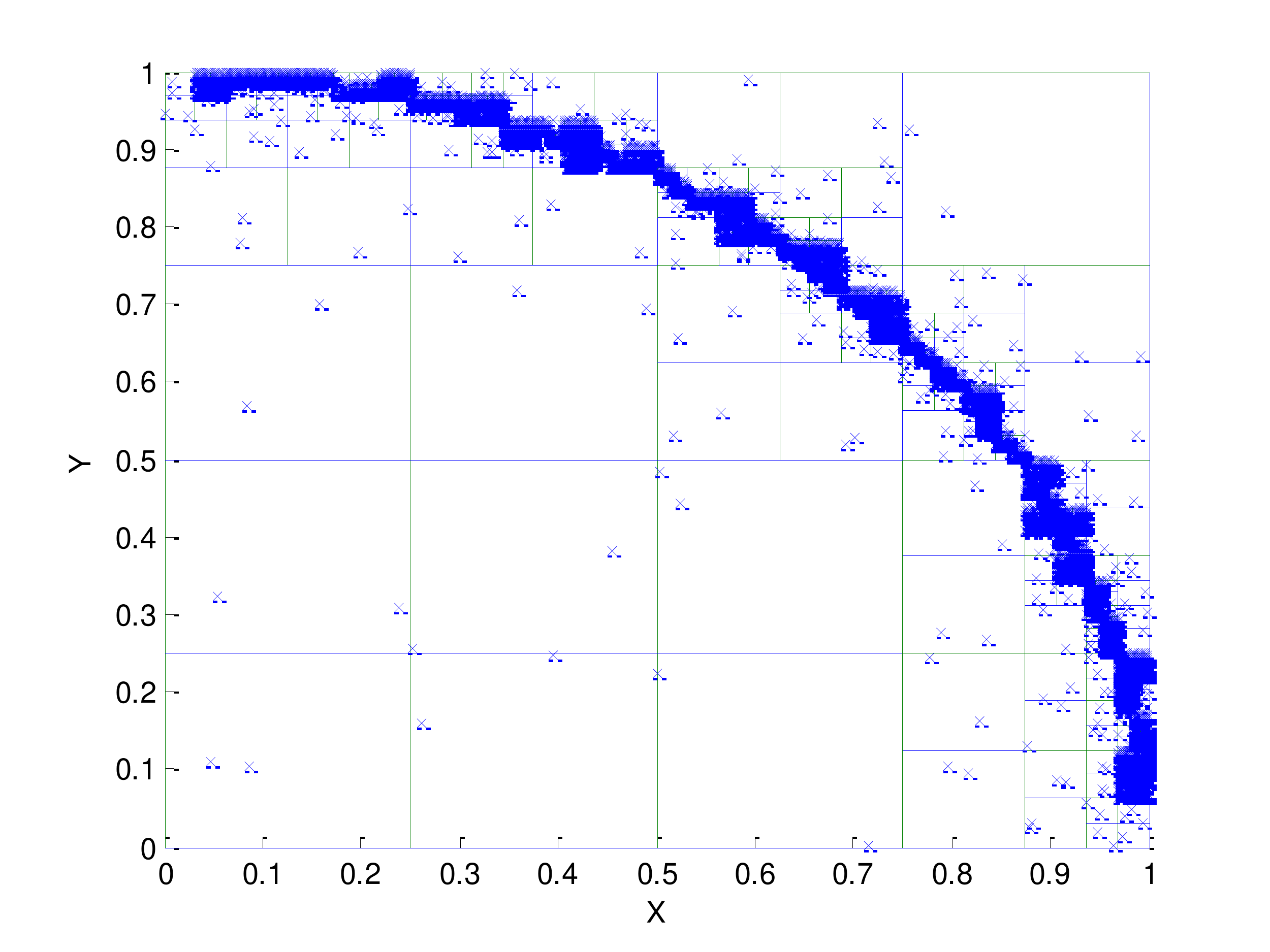}\vspace{-.7cm}
  \caption{ repartition of the points  (sixth iteration).} \label{step8}
 \end{minipage}
 \end{figure}

\noindent We begin the numerical tests with the algorithm (Alog
1). Figures \ref{step1}-\ref{step8} show for $N=10000$ and $L=6$
the evolution of the mesh and the repartition of the random points
during the iterations.  We remark that this points are
concentrated around the curve $x^2+y^2=1$ where the function $f_c$
represents a singularity.

\begin{figure}[h!]
%\vspace{.5cm}
 \begin{minipage}[h]{.46\linewidth}\vspace{-.9cm}
  \centering\includegraphics[width=7.5cm]{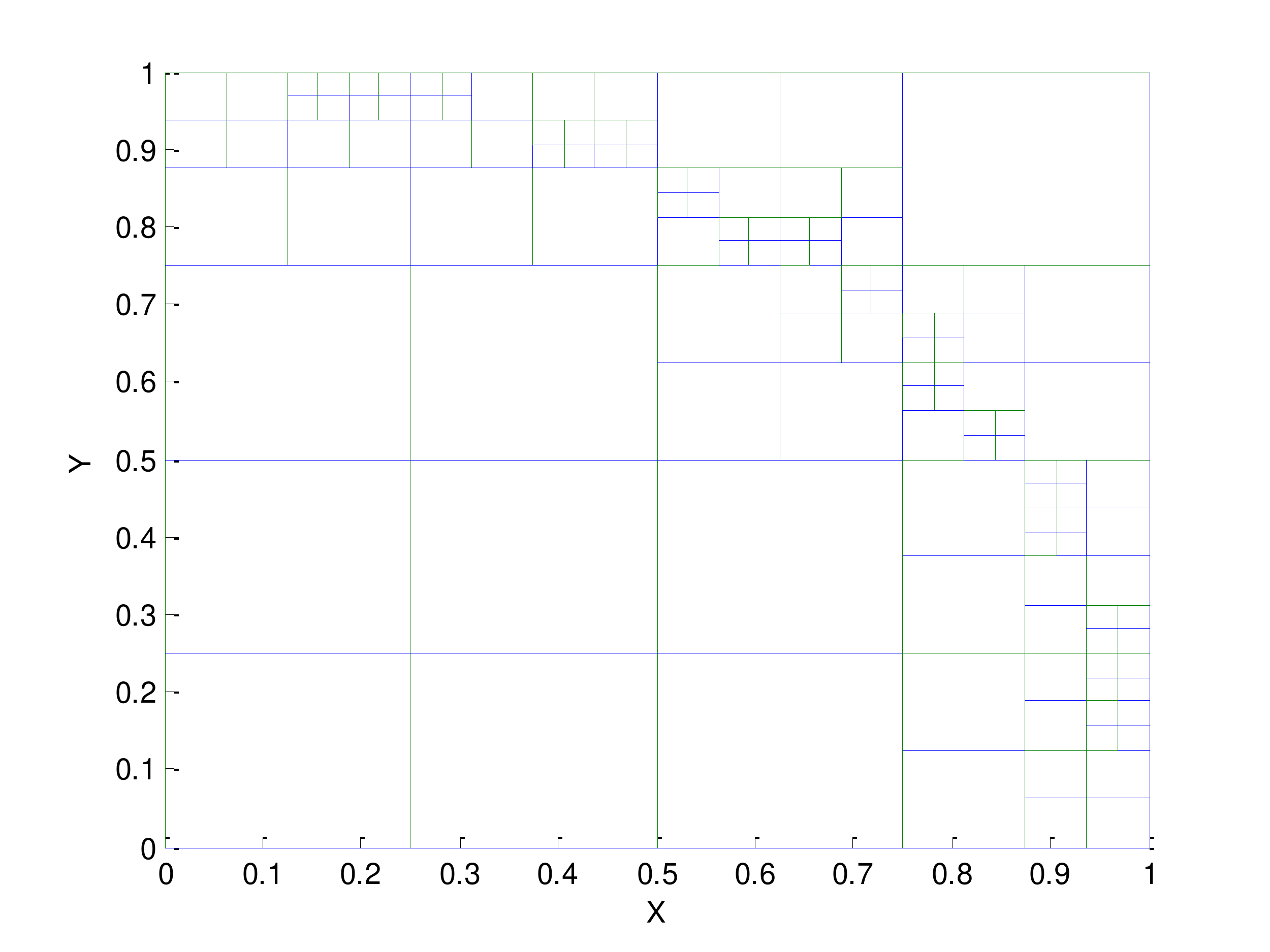}\vspace{-.7cm}
  \caption{Mesh for the fourth iteration.}
 \end{minipage} \hfill
 \begin{minipage}[h]{.46\linewidth}\vspace{-.6cm}
  \centering\includegraphics[width=7.5cm]{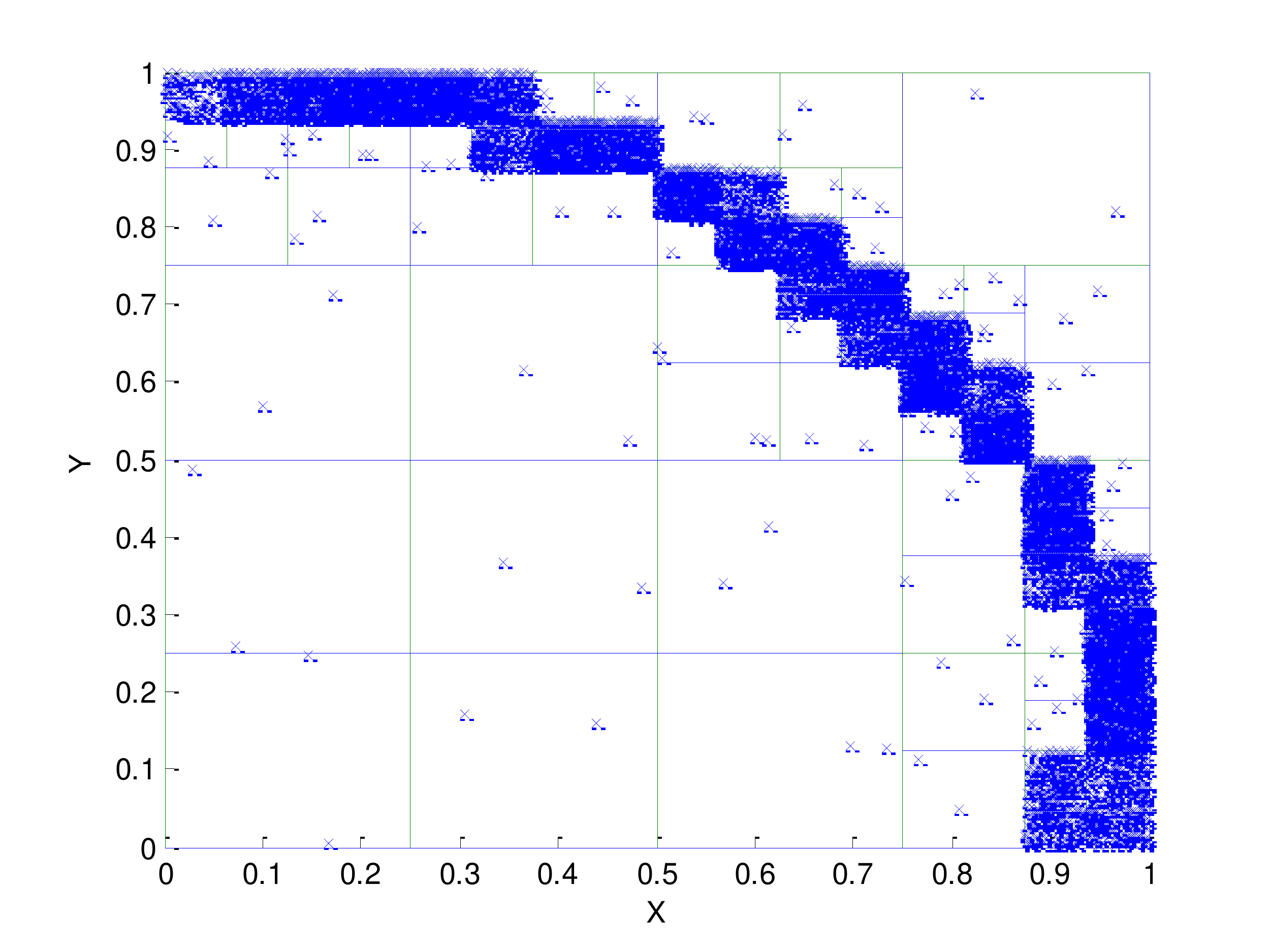}\vspace{-.7cm}
  \caption{ repartition of the points (fourth iteration).}
 \end{minipage}
\end{figure}

\begin{figure}[h!]
%\vspace{.5cm}
 \begin{minipage}[h]{.46\linewidth}\vspace{-.4cm}
  \centering\includegraphics[height=5.8cm,width=8.3cm]{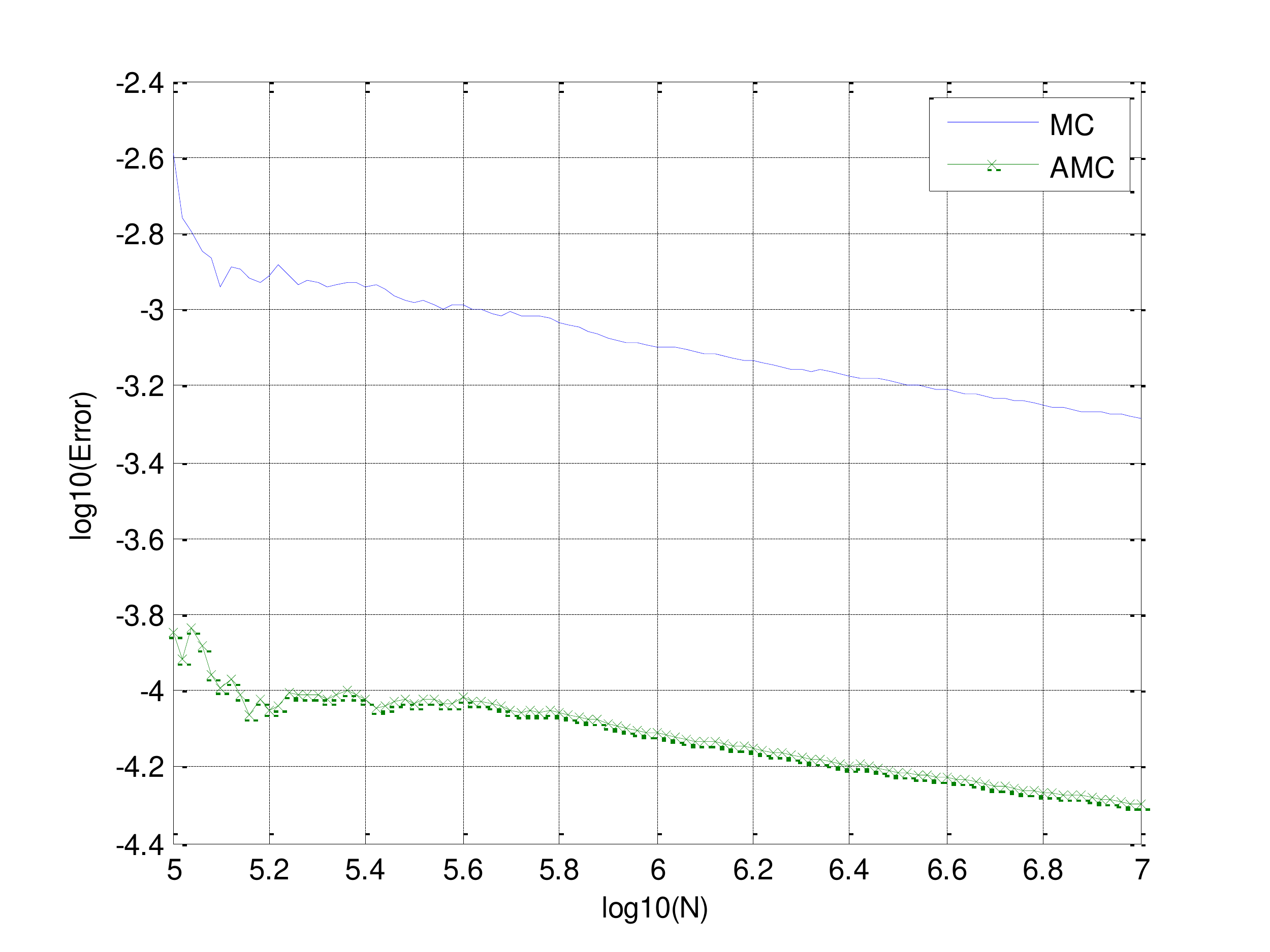}\vspace{-.5cm}
  \caption{First test case: $E_{MC}$ and $E_{AMC}$ with respect to $N$ in logarithmic scale.} \label{com_vol_exact_num_disc}
 \end{minipage} \hfill
 \begin{minipage}[h]{.46\linewidth}\vspace{-.0cm}
  \centering\includegraphics[height=5.8cm,width=8.3cm]{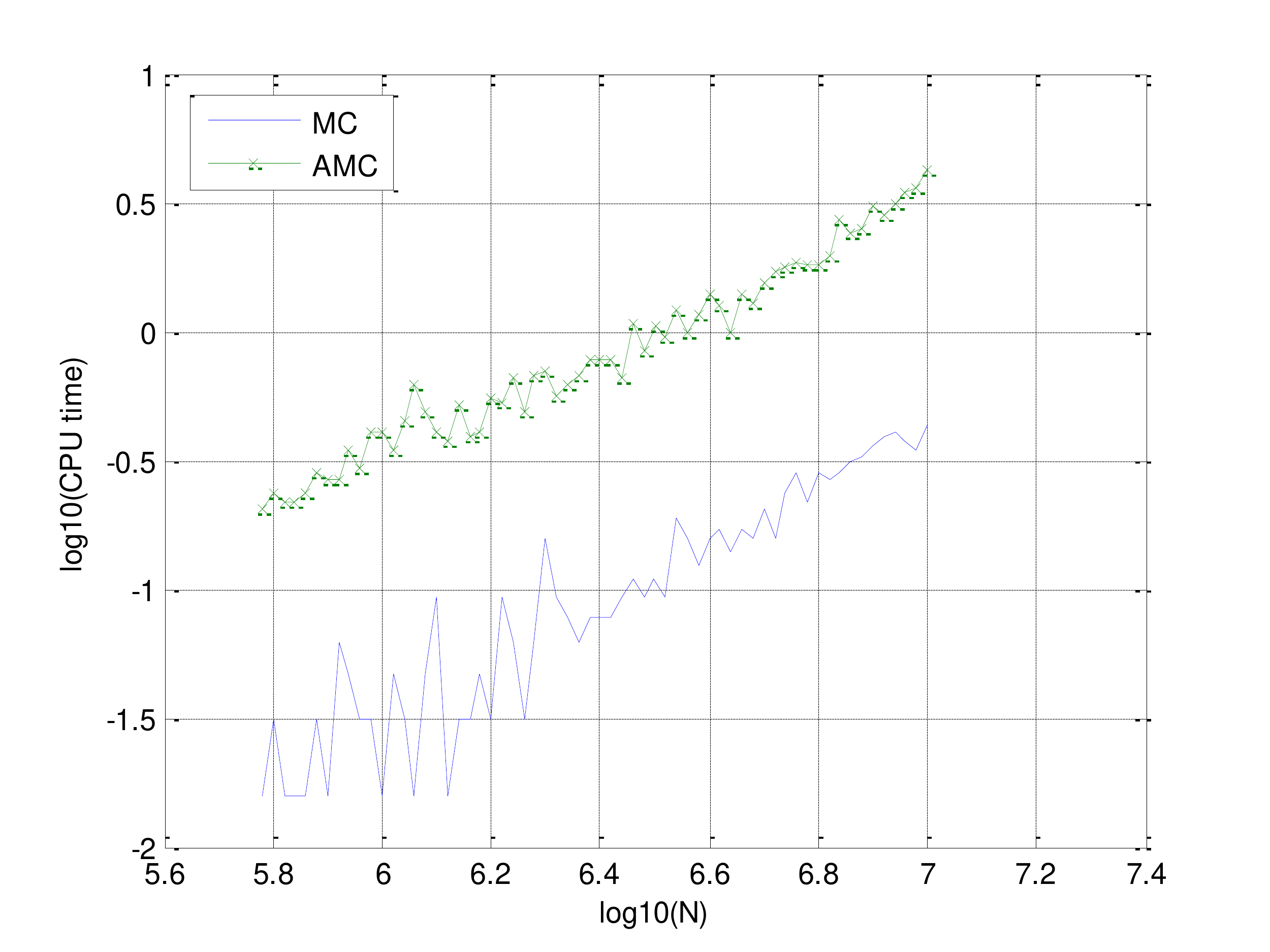}\vspace{-.5cm}
  \caption{First test case: CPU time of the MC and AMC methods with respect to $N$ in logarithmic scale.} \label{com_time_exact_num_disc}
 \end{minipage}
\end{figure}

Figure \ref{com_vol_exact_num_disc} shows a comparison in
logarithmic scale of the relative errors ($I=\ds \frac{\pi}{4}$)
$$ E_{MC} = \ds \frac{I - I_{MC}}{I}$$
corresponding to the MC method and
$$ E_{AMC} = \ds \frac{I - I_{AMC}}{I}$$
corresponding to the AMC method with respect to the number of
random points $N$ where the total number of the iteration $L=4$.
As we can see in figure \ref{com_vol_exact_num_disc}, the AMC
method is more precise than the MC method. Still we have to
compare the efficiency of the AMC method with respect to the CPU
time of computation. In fact, figure \ref{com_time_exact_num_disc}
shows that for the considered $N$, the corresponding CPU times for
the AMC are smaller from those with MC. In particular, the MC
method gives for $N=10^7$ an error of $E_{MC}=0.00052$ with a CPU
time of $0.44 s$, but the AMC gives for $N=10^6$ an error of
$E_{AMC}=0.00008$ with a CPU time of $0.4 s$. Hence, the powerful
of the AMC method. It is also clear that to get more precision
with the AMC method, we can increase the number of iterations $L$.

Next, we consider the algorithm (Algo 2) with $N_{ess}=100$,
$L=4$.

\begin{figure}[h!]
%\vspace{.5cm}
 \begin{minipage}[h]{.46\linewidth}\vspace{-.2cm}
  \centering\includegraphics[height=5.8cm,width=8.3cm]{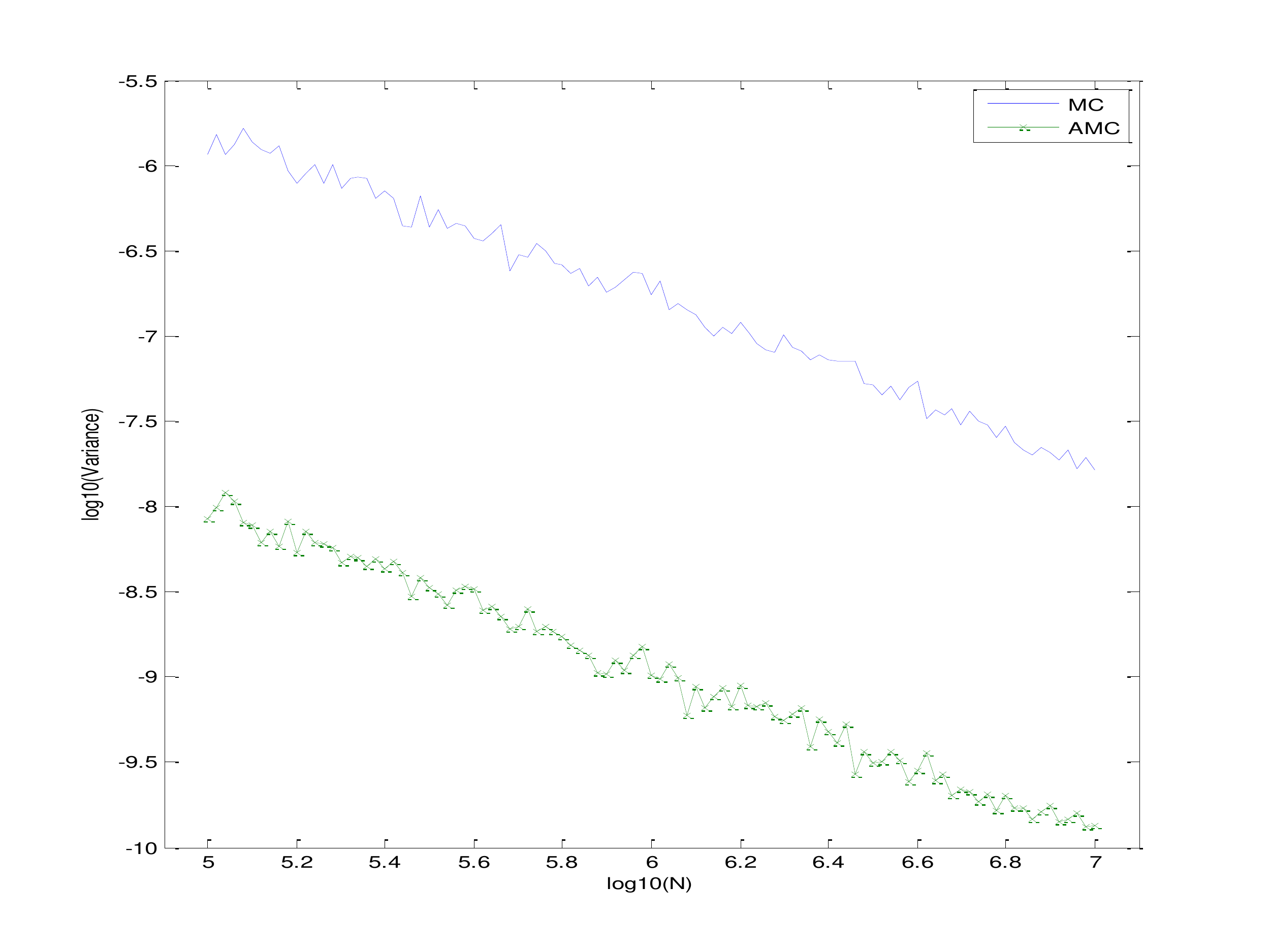}\vspace{-.5cm}
  \caption{First test case: Estimated variances $V_{MC}$ and $V_{AMC}$ with respect to $N$ in logarithmic scale.} \label{variance_cercle}
 \end{minipage} \hfill
 \begin{minipage}[h]{.46\linewidth}\vspace{-.2cm}
  \centering\includegraphics[height=5.8cm,width=8.3cm]{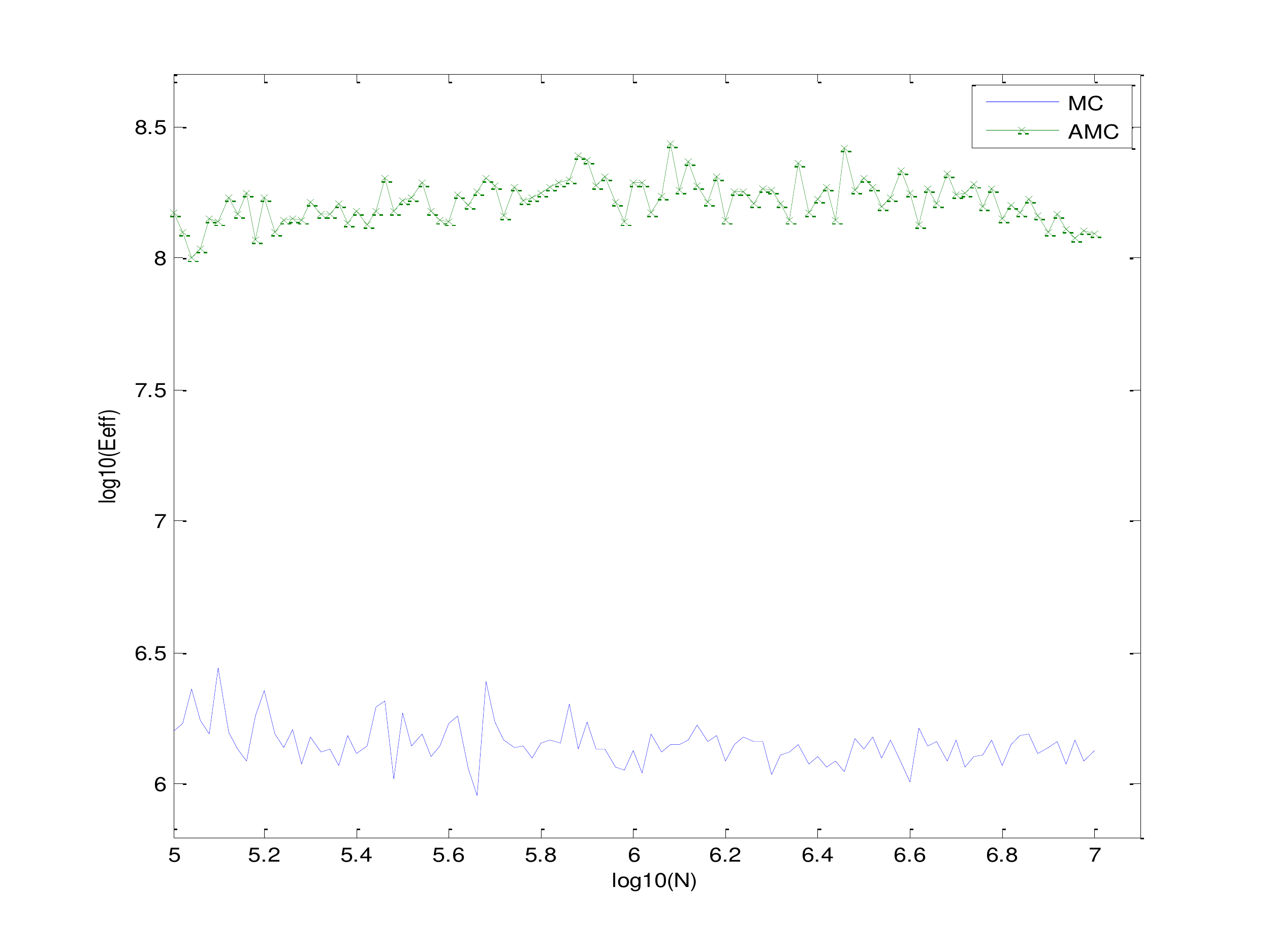}\vspace{-.5cm}
  \caption{First test case: Efficiencies $E^{eff}_{MC}$ and $E^{eff}_{AMC}$ with respect to $N$ in logarithmic scale.} \label{efficacite_cercle}
 \end{minipage}
\end{figure}

\noindent Figure \ref{variance_cercle} shows the comparison of the
estimated  variance between the classical Monte Carlo ($V_{MC}$)
and adaptive Monte Carlo method ($V_{AMC}$) in logarithmic scale.
As the adaptive algorithm consists to minimize the variance, it is
clear in this figure that the goal is attended.  Figure
\ref{efficacite_cercle} shows in logarithmic scale the efficiency
of the MC and AMC methods versus the number of random points $N$
by using the following formulas (see \cite{PE})
\[
E^{eff}_{MC} = \ds \frac{1}{T_{MC}*V_{MC}}
\]
and
\[
E^{eff}_{AMC} = \ds \frac{1}{T_{AMC}*V_{AMC}},
\]
where $T_{MC}$ and $T_{AMC}$ are respectively the CPU time of the
MC and AMC methods. It is clear that the efficiency of the AMC
method is more important than the MC method.
\subsubsection{Second test case} 
\begin{figure}[h!]

 \begin{minipage}[h]{.46\linewidth}\vspace{-.3cm}
  \centering\includegraphics[width=7.cm]{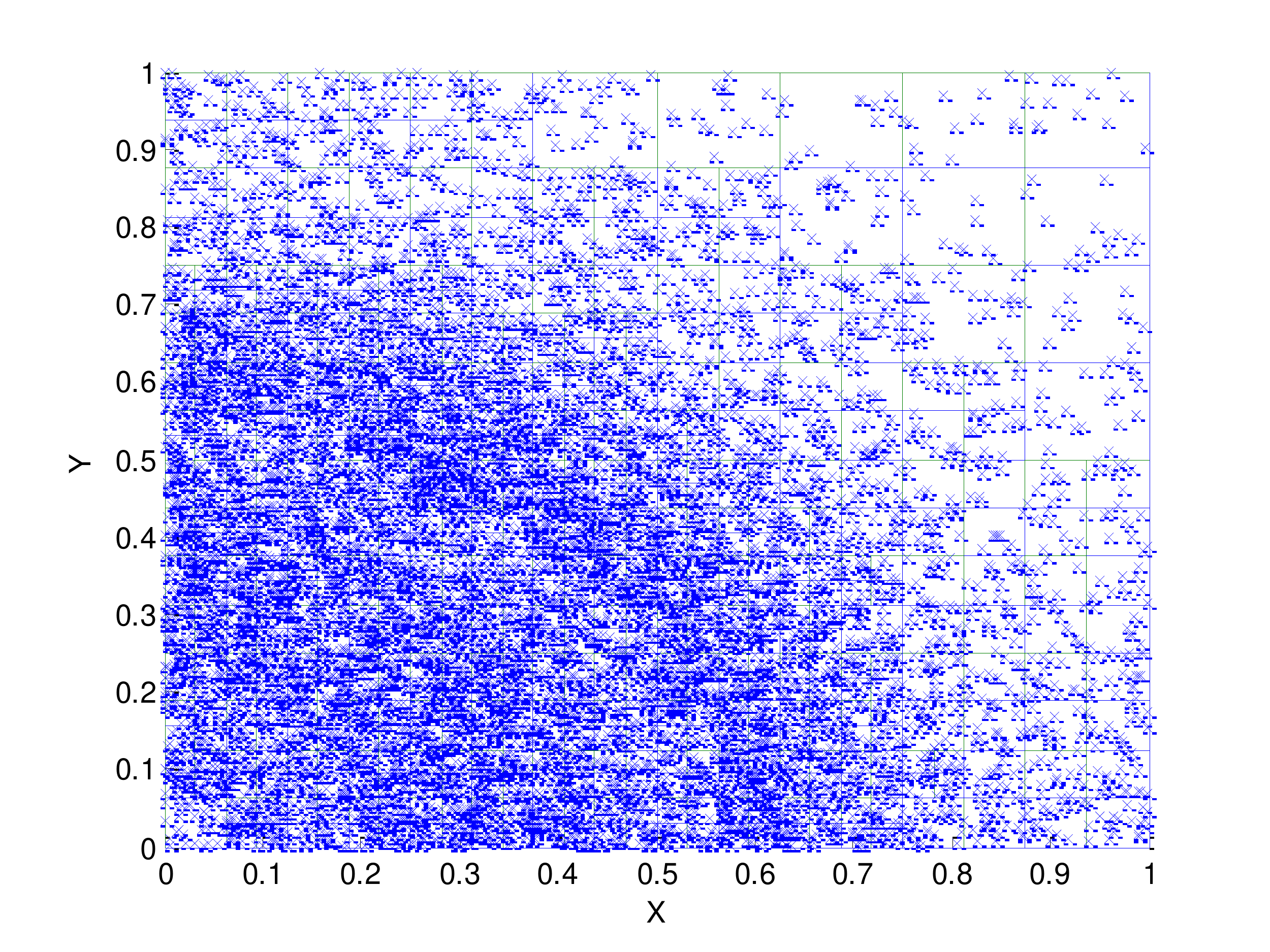}\vspace{-.5cm}
  \caption{AMC mesh, $\alpha=-5$.}  \label{alpham5}
 \end{minipage} \hfill
 \begin{minipage}[h]{.46\linewidth}\vspace{-.3cm}
  \centering\includegraphics[width=7.cm]{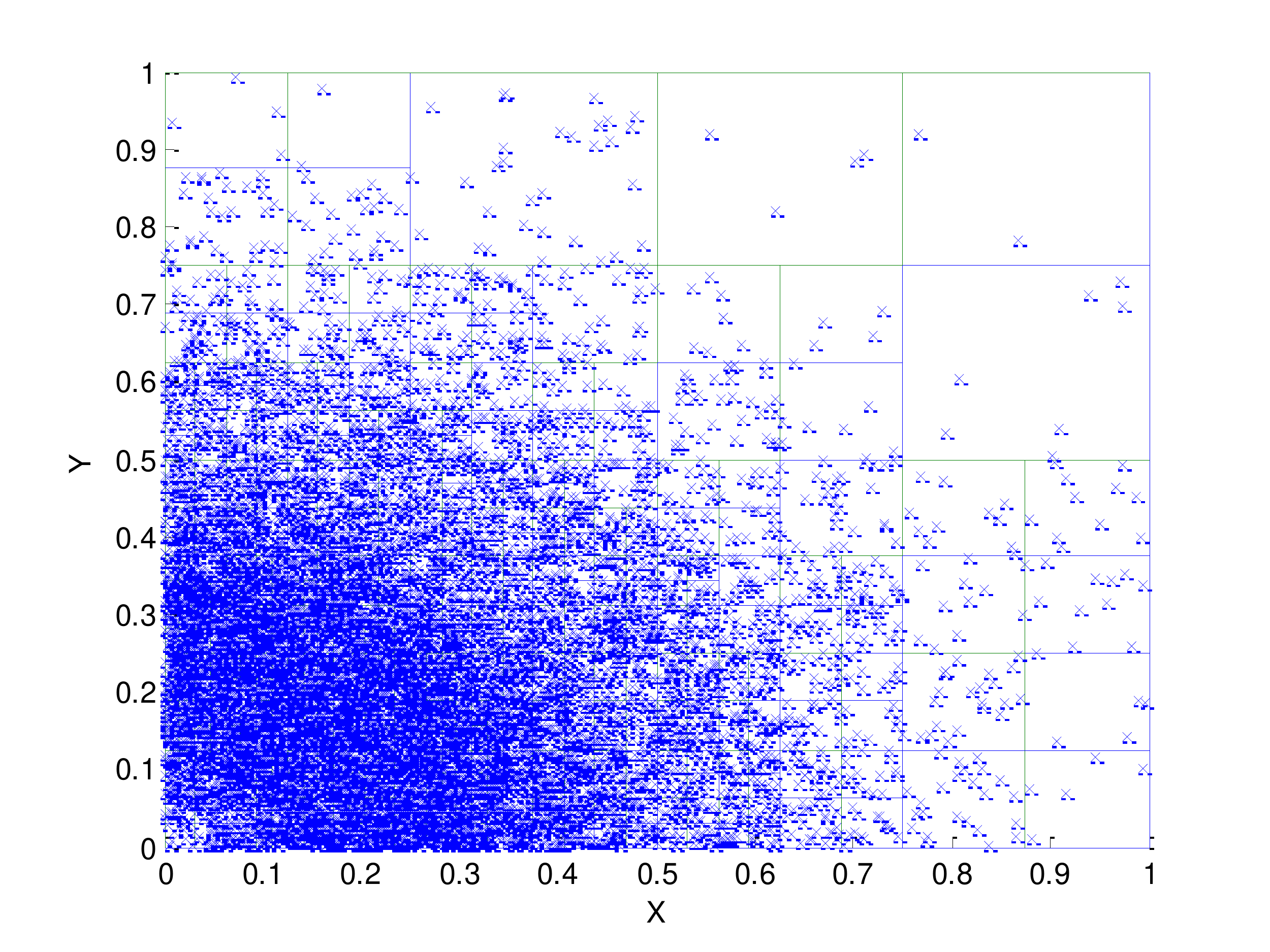}\vspace{-.5cm}
  \caption{AMC mesh $\alpha=-10$.}
 \end{minipage}

\vspace{.3cm}
 \begin{minipage}[h]{.46\linewidth}\vspace{-.2cm}
  \centering\includegraphics[width=7.cm]{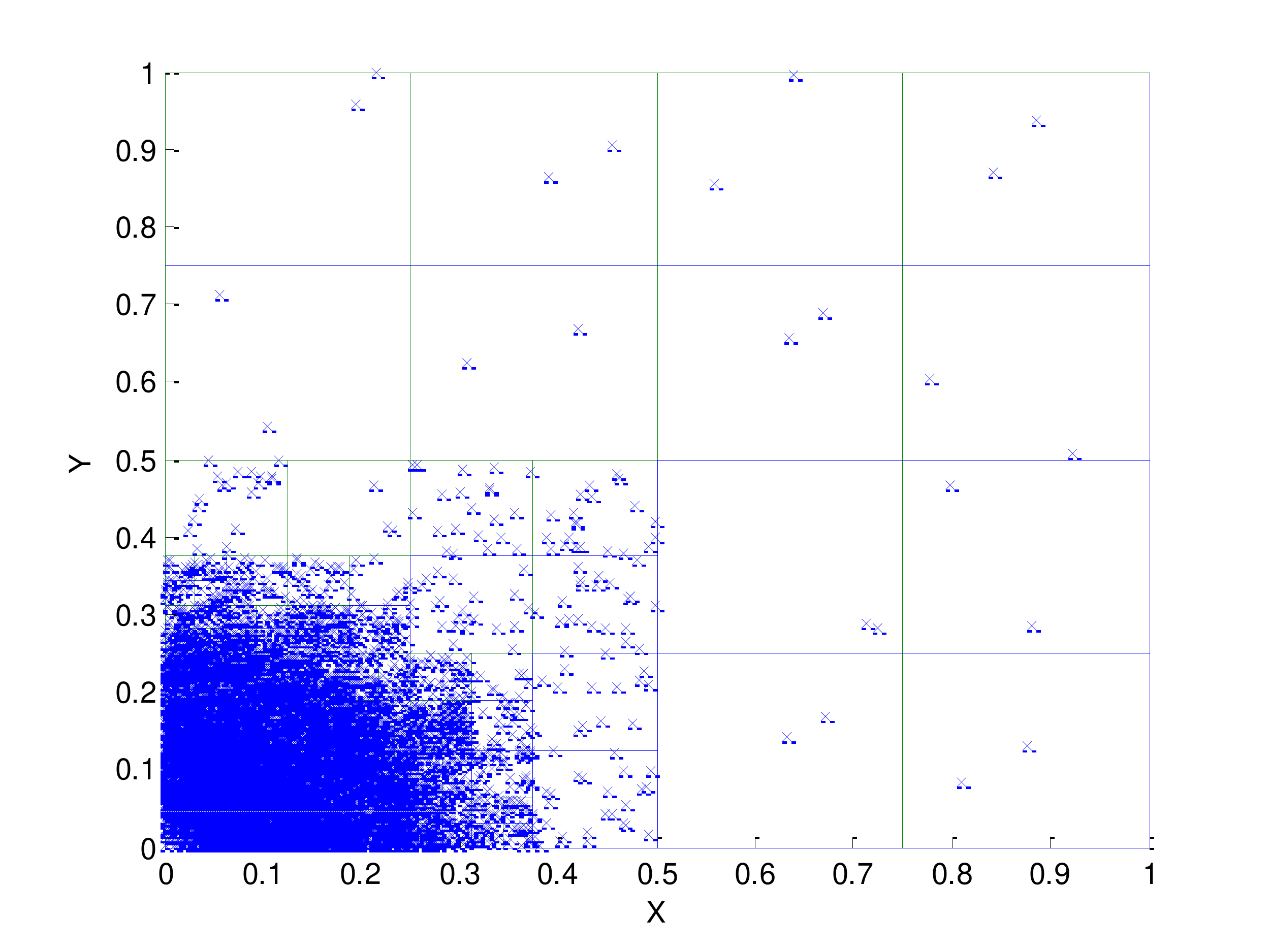}\vspace{-.5cm}
  \caption{AMC mesh $\alpha=-50$.}
 \end{minipage} \hfill
 \begin{minipage}[h]{.46\linewidth}\vspace{-.2cm}
  \centering\includegraphics[width=7.cm]{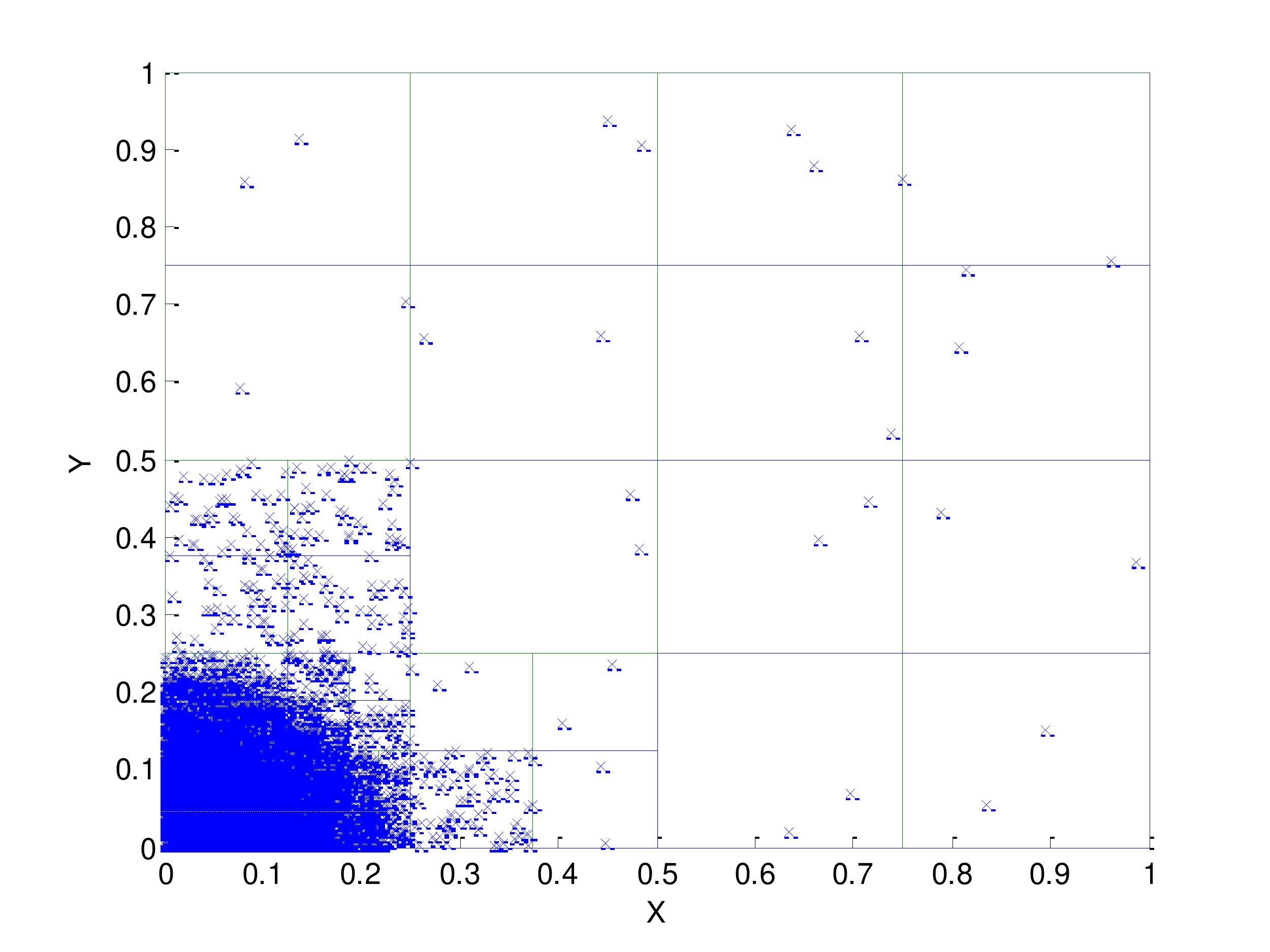}\vspace{-.5cm}
  \caption{AMC mesh $\alpha=-100$.} \label{alpham100}
 \end{minipage}

 \end{figure}

$\mbox{ }$\\
\noindent In this case, we consider the function
\[
f_{2,g}(x,y) = e^{-\alpha (x^2+y^2)},
\]
where $\alpha$ is a real positive parameter. We begin the adaptive
algorithm with the same initial mesh as the previous case and we
choose $N=10000$. Figures \ref{alpham5}-\ref{alpham100} show for
$L=6$ the meshes and random points repartition with respect to
$\alpha$. When $\alpha$ increase, the mesh and the random points
follow the function $f_{2,g}$ and focus more and more around the
origin
of axis.\\

\newpage

\begin{figure}[h!]
%\vspace{.5cm}
 \begin{minipage}[h]{.46\linewidth}\vspace{-.4cm}
  \centering\includegraphics[height=5.8cm,width=8.3cm]{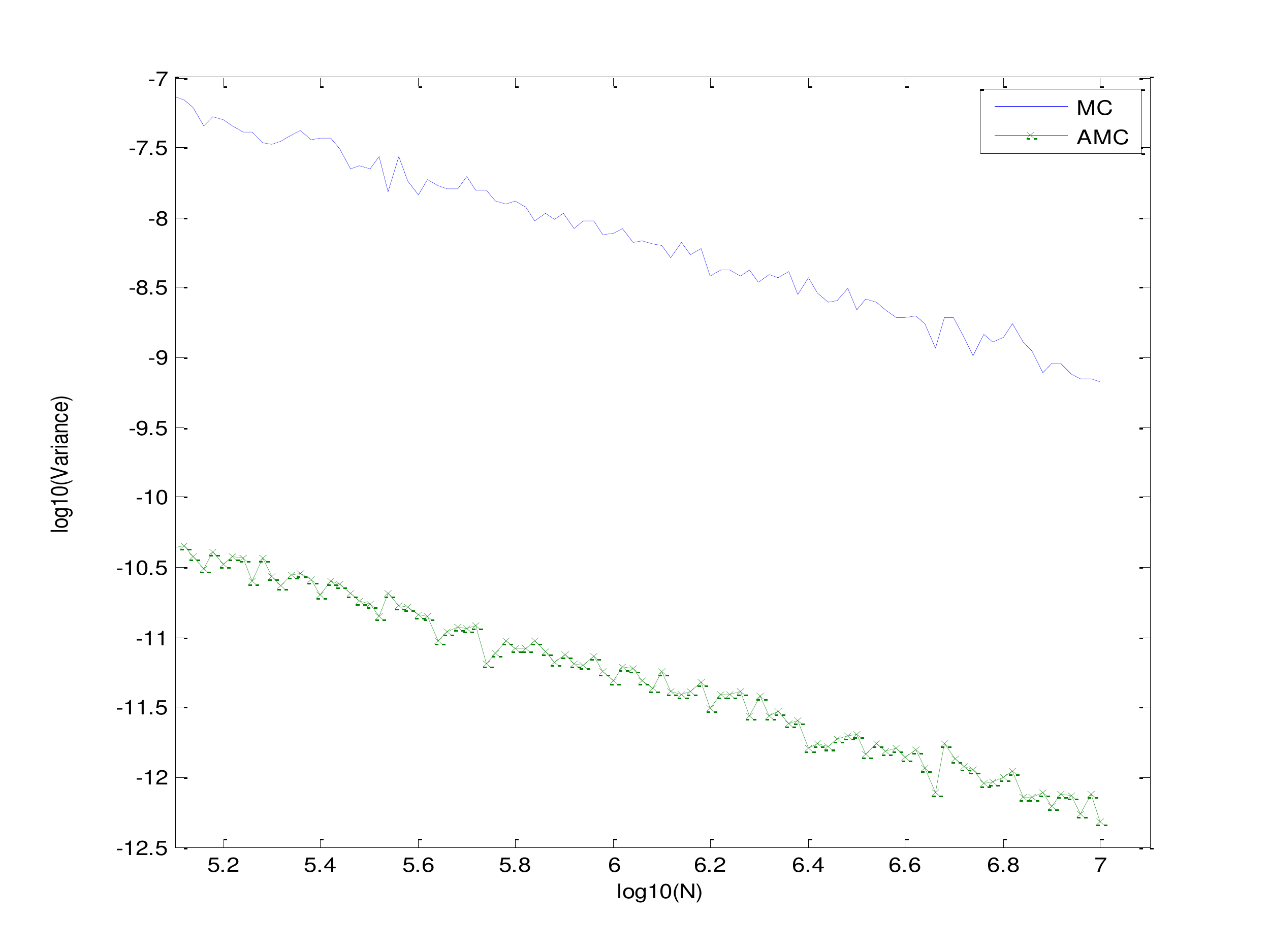}\vspace{-.5cm}
  \caption{Second test case ($\alpha=-50$): Estimated variances $V_{MC}$ and $V_{AMC}$ with respect to $N$ in logarithmic scale.} \label{variance_exp}
 \end{minipage} \hfill
 \begin{minipage}[h]{.46\linewidth}\vspace{-.4cm}
  \centering\includegraphics[height=5.8cm,width=8.3cm]{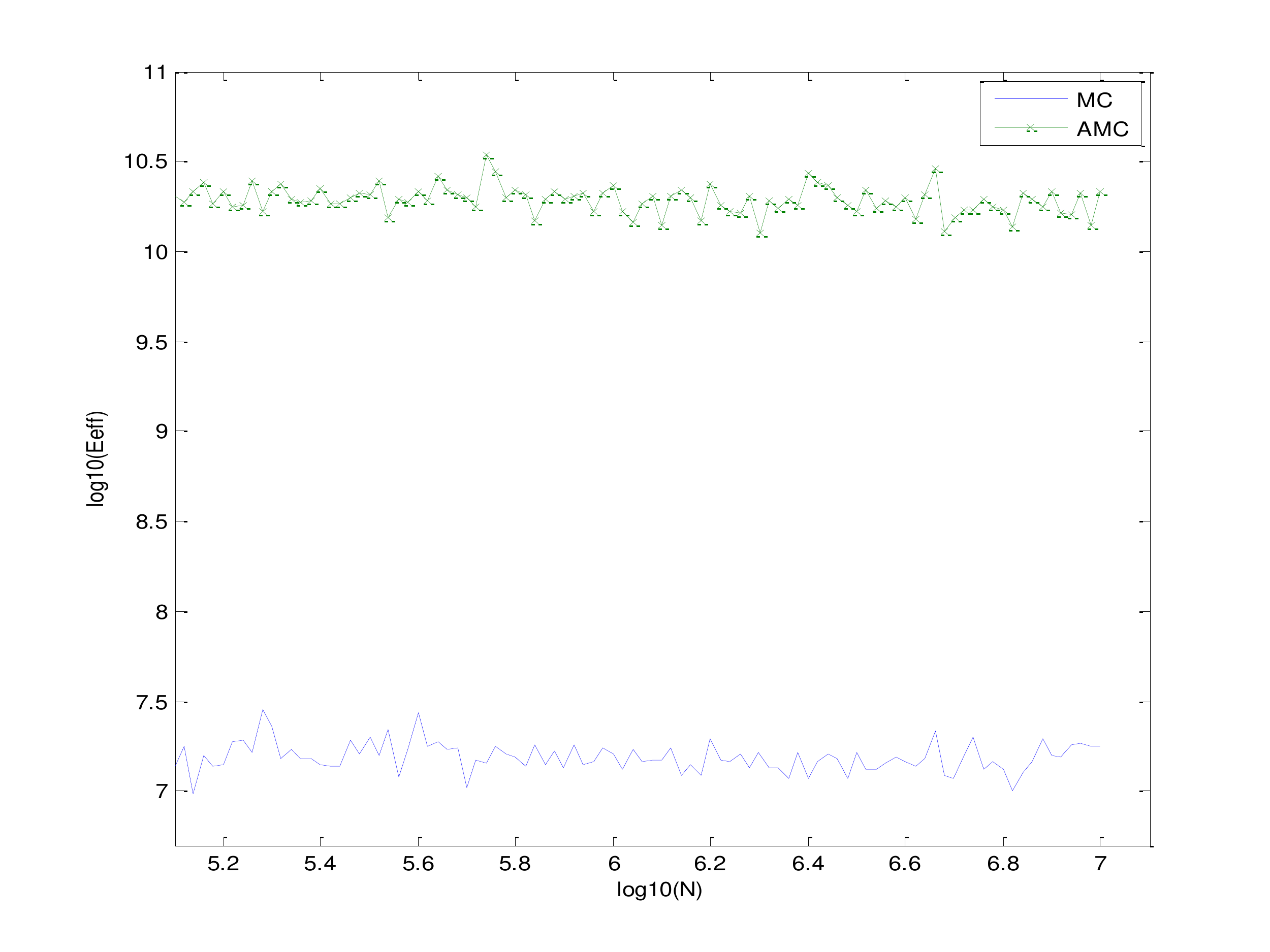}\vspace{-.5cm}
  \caption{Second test case ($\alpha=-50$): Efficiencies $E^{eff}_{MC}$ and $E^{eff}_{AMC}$ with respect to $N$ in logarithmic scale.} \label{efficacite_exp}
 \end{minipage}
\end{figure}

\noindent Figure \ref{variance_exp} shows for $\alpha=-50$,
$N_{ess}=100$ and $L=4$, the comparison of the estimated  variance
between MC and AMC methods with respect to $N$ in logarithmic
scale. Figure \ref{efficacite_exp} shows in logarithmic scale the
efficiency of the MC and AMC methods versus the number of points
$N$. One more time, it is clear that the efficiency of the AMC
method is more important than the MC one.
\newpage
\subsection{3D validations}
\begin{figure}[h!]
\includegraphics[width=8cm]{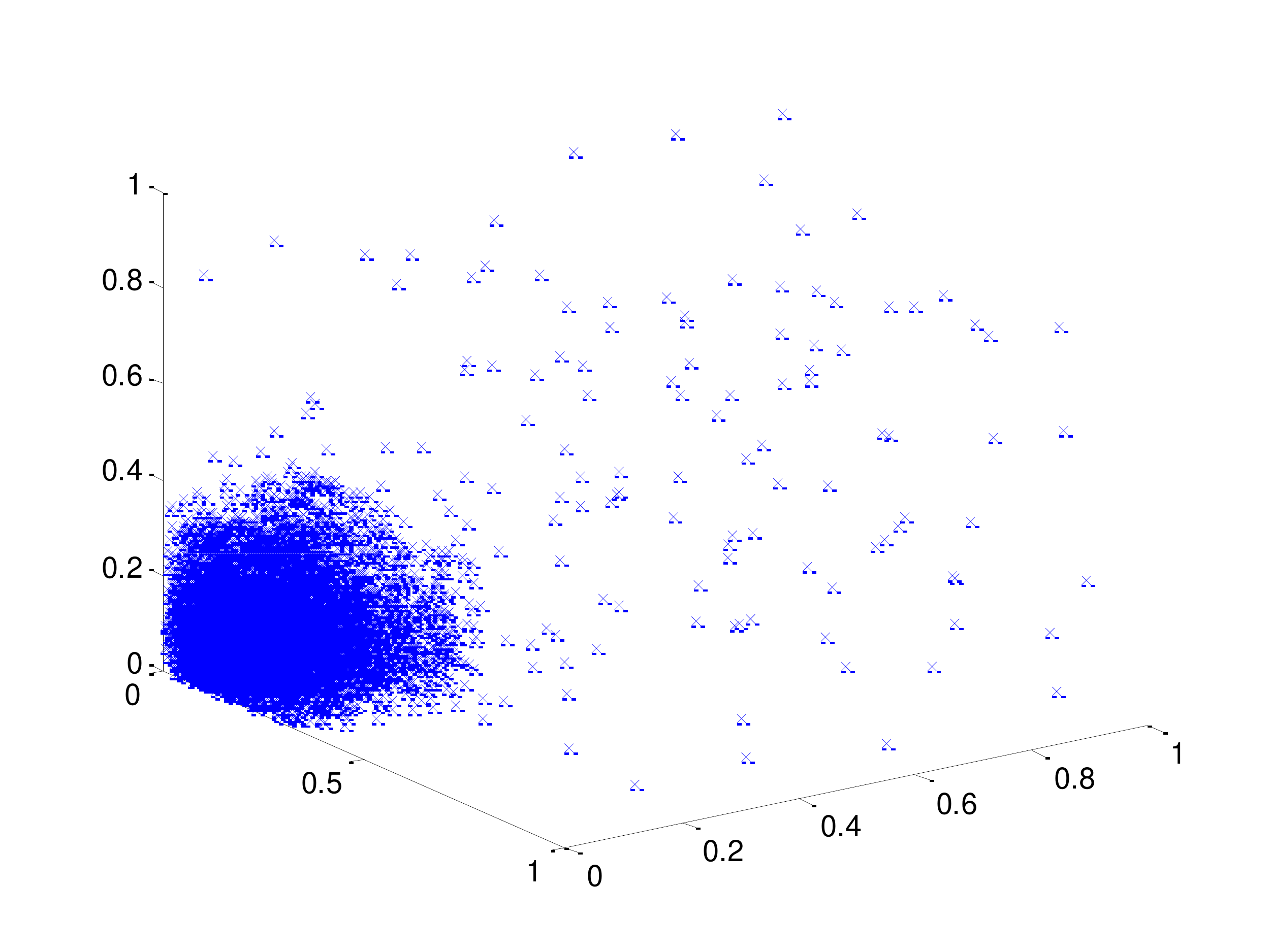}\vspace{-.5cm}
\caption{Mesh Gaussian for $\alpha=-50$.}  \label{mg3D}
\end{figure}

In this section, we consider the unit cube $D=[0,1[^3$. We
consider the function
\[
f_{3,g}(x,y) = e^{-\alpha (x^2+y^2+z^2)},
\]
where $\alpha$ is a real positive parameter. The initial mesh is
constituted by a regular partition with $N_0=4$ segments in every
side of $D^1=D$. Figure \ref{mg3D}
shows the repartition of the random points for $\alpha=-50$, $L=6$ and $N=10000$.\\

\begin{figure}[h!]
%\vspace{.5cm}
 \begin{minipage}[h]{.46\linewidth}\vspace{-.4cm}
  \centering\includegraphics[height=5.8cm,width=8.3cm]{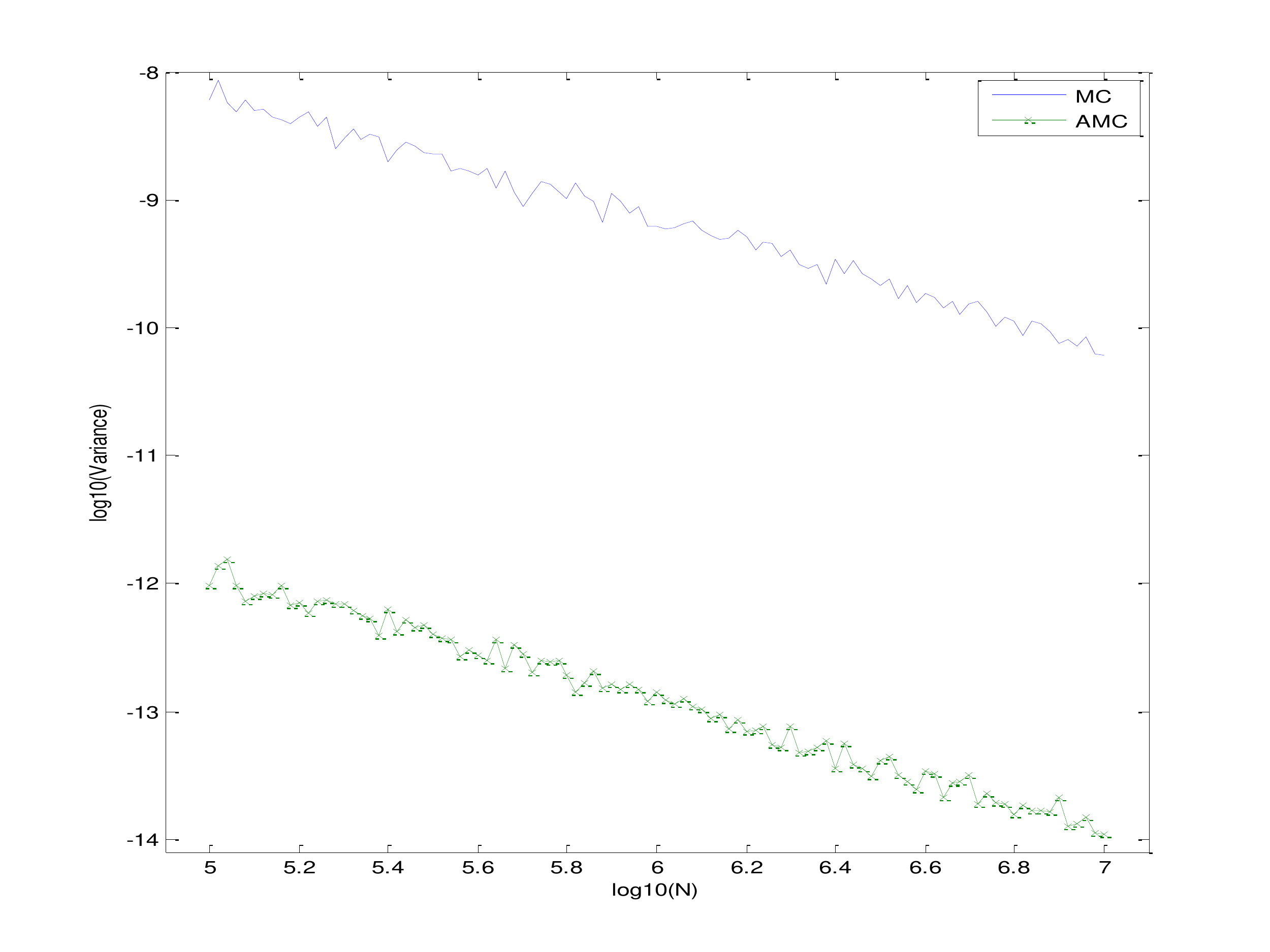}\vspace{-.5cm}
  \caption{3D test case ($\alpha=-50$): Estimated variances $V_{MC}$ and $V_{AMC}$ with respect to $N$ in logarithmic scale.} \label{variance3D}
 \end{minipage} \hfill
 \begin{minipage}[h]{.46\linewidth}\vspace{-.4cm}
  \centering\includegraphics[height=5.8cm,width=8.3cm]{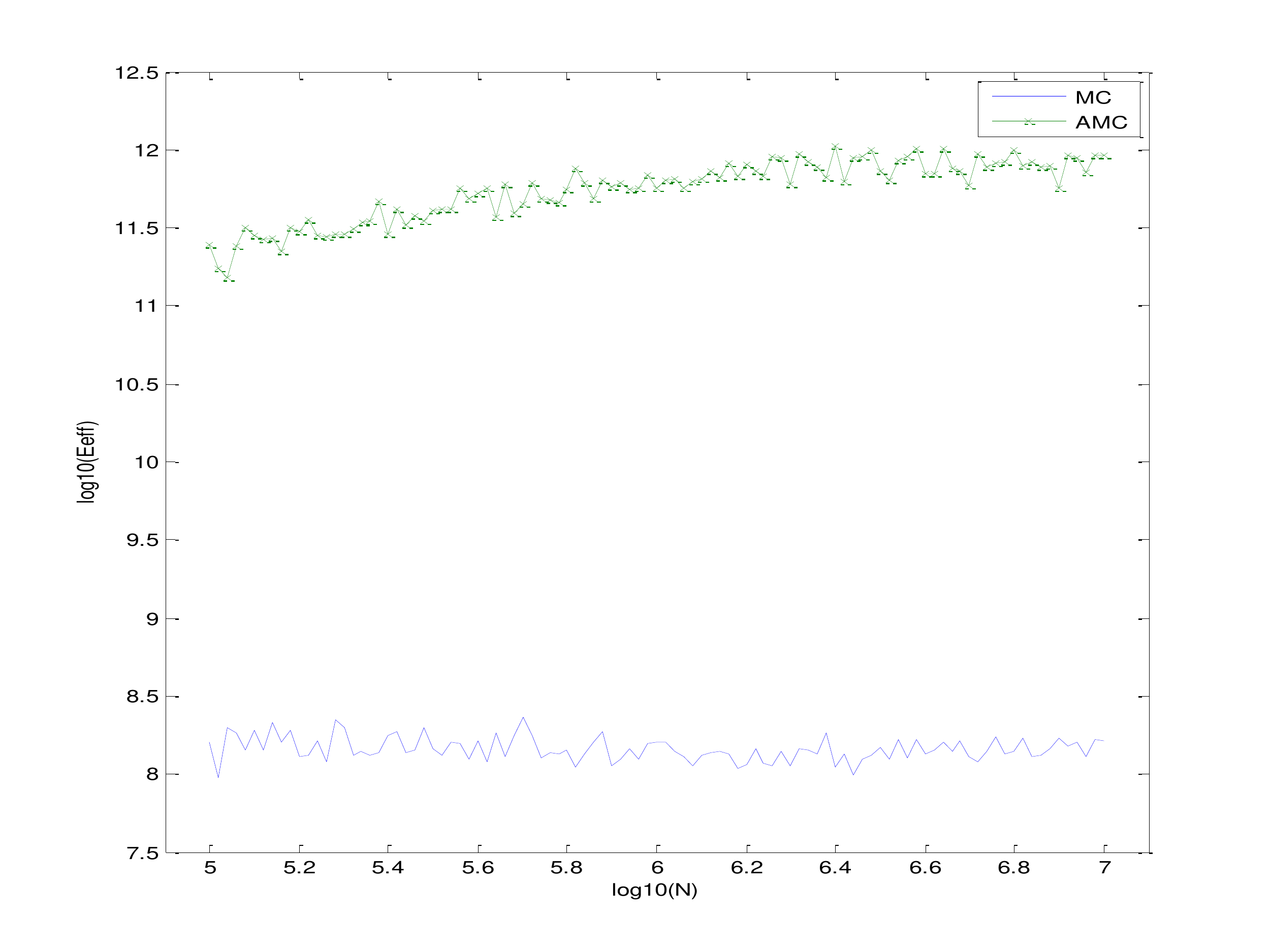}\vspace{-.5cm}
  \caption{3D test case ($\alpha=-50$): Efficiencies $E^{eff}_{MC}$ and $E^{eff}_{AMC}$ with respect to $N$ in logarithmic scale.} \label{efficacite_gauu_3D_exp}
 \end{minipage}
\end{figure}

\noindent As for the previous case, figure \ref{variance3D} shows
for $N_{ess}=100$ and $L=4$, the comparison of the estimated
variance between MC and AMC methods with respect to $N$ in
logarithmic scale. Figure \ref{efficacite_gauu_3D_exp} shows in
logarithmic scale the efficiency of the MC and AMC methods versus
the number of points $N$. We can deduce the same remark for the
efficiency of the AMC method in dimension three.

{{\bf{Acknowledgments}}  We are grateful to my colleague Rami El
Haddad for all the help given by him.}
\newpage

\end{document}